\documentclass[a4paper,twocolumn]{article} 
                          
\usepackage{graphicx}   
\usepackage{epsfig} 
\usepackage{epstopdf}
\usepackage{times} 
\usepackage{amsmath} 
\usepackage{amssymb}  
\usepackage{theorem}
\usepackage[]{units}
\usepackage{subfig}
\usepackage{psfrag}
\usepackage{wrapfig}
\usepackage{color}
\usepackage[percent]{overpic}

\usepackage[dvipsnames]{xcolor}

\usepackage[]{algorithm}
\usepackage[]{algpseudocode}

\usepackage{multirow}

\newcommand{\R}{\mathbb{R}}

\newtheorem{lemma}{Lemma}
\newtheorem{proof}{Proof}
\newtheorem{proposition}{Proposition}

\newtheorem{example}{Example}

\usepackage[font={footnotesize}]{caption} 

\title{
Regional predictive control with suboptimally extended regions of validity
}

\author{Kai K\"onig and Martin M\"onnigmann\thanks{Corresponding author.}\\
Automatic Control and Systems Theory, Department of Mechanical Engineering,\\
	Ruhr-Universit\"at Bochum, 44801 Bochum, Germany.\\ E-mail: {\tt\small kai.koenig-h4d@rub.de} and {\tt\small martin.moennigmann@rub.de}}

\begin{document}

\maketitle

\begin{abstract}  
Model predictive control (MPC) is based on perpetually solving optimization problems. 
The solution of the optimization is usually interpreted as the optimal input for the current state. However, the solution of the optimization does not just provide an optimal input, but an entire optimal affine feedback law and a polytope on which this law is optimal. We recently proposed to use this feedback law as long as the system remains in its polytope. This can be interpreted as an event-based approach, where leaving the current polytope is the event that triggers the next optimization. 
This approach is especially appropriate for a networked control setting since the feedback laws and their polytopes can be evaluated with a low computational effort on lean local nodes. In this article the region of validity for a feedback law is extended. More precisely, the optimal polytopes are extended to nonlinearly bounded regions of validity resulting from the intersection of a feasibility and stability region. As a result, fewer quadratic programs need to be solved compared to the optimal approach. The new validity regions are still suitable for the evaluation on a lean local node. Moreover, the regions can be adjusted for a desired closed-loop performance.
\end{abstract}

\section{Introduction}
Model predictive control is based on the solution of an optimal control problem (OCP) in every time step. Many approaches have been presented for overcoming the high computational effort resulting from the online optimizations. While one group focuses on speeding up the numerical solution of the OCP (see, e.g., \cite{Ferreau2008, Pannocchia2007}) another group proposes to avoid online optimizations whenever possible (see, e.g., \cite{Bemporad2002, Kvasnica2017, Gupta2011, Oberdieck2017, Tondel2001}). 
Event-triggered MPC belongs to the second group and proposes to apply feedback not periodically but only when the system needs attention. The event that triggers a new optimization can, e.g., be the difference between the predicted and actual trajectory as proposed in \cite{Bernardini2012, Ferrara2014, Incremona2015, Li2014, Lehmann2013} or the rate of change of the cost function as presented in \cite{Eqtami2010, Eqtami2011}.  \\
Regional predictive control is an event-triggered MPC approach that exploits the structure of the solution to the linear MPC problem as triggering event (see \cite{Jost2015a, Koenig2017d, Koenig2017c, Koenig2017a, MonnigmannOttenNMPC2015, SchulzeDarup2017CDC}). It is based on the fact that the solution at a point provides an affine piece, i.e., a feedback law and its polytopic region of validity, with low effort. Instead of optimizing in every time step, feedback laws are calculated and reused as long as the closed-loop system remains in the corresponding polytope. 
Leaving the current polytope is the event that triggers the solution of a new quadratic program (QP). The regional MPC approach is appropriate for a networked control setting, where solving  QPs can be outsourced to a powerful central node. The affine piece can be transmitted to a local node where it can be evaluated with a much smaller computational effort than required for the QP, thus requiring only a lean local node (see \cite{BernerP2016}). \\
Obviously, it is of special interest to extend the region of validity for a feedback law to reduce the number of solved QPs further. By dropping the requirement of optimality, a feedback law can be used even outside its original polytope as long as it is feasible and stabilizing. The idea was already proposed in \cite{Jost2015a}, but the approach presented there was not suitable for a networked control setting with a lean local node, since it required a computationally expensive evaluation of feasibility and stability criteria. 
Results on how to extend the suboptimal approach from \cite{Jost2015a} for use in a networked setting with lean local nodes were first presented in \cite{Koenig2017c}. The present article extends the approach from \cite{Koenig2017c} by an alternative calculation of the feasibility region, and by a closed-loop performance analysis.   \\
We determine extended regions of validity for a feedback law that are suitable for the use in a networked control setting. It is the idea to calculate an explicit nonlinearly bounded state-space region on which a feedback law is feasible and stabilizing. Such a region results from the intersection of a polytope  and another region described by a simple quadratic inequality. 
The new regions reduce the number of QPs to be solved compared to the optimal approach \cite{Jost2015a} and, despite being nonlinearly bounded, it is computationally cheap to evaluate them, and thus they can be evaluated on a low-cost node. Furthermore, the new regions can be adjusted with respect to closed-loop performance.\\
We state the system and problem class along with some preliminaries in Section~\ref{sec:problemstatement}. The new approach is presented in Section~\ref{sec:mainpart} and applied to two examples in Section~\ref{sec:example}. A short summary is given in Section~\ref{sec:conclusion}. 

\section{Problem statement and preliminaries}\label{sec:problemstatement}
We consider the optimal control problem
\begin{subequations}\label{eq:MPCProblem}
\begin{align}
&\min \limits_{\substack{u(k),k=0,\ldots,N-1\\x(k),k=1,\ldots,N}}  \quad  \| {x(k)}\|_P^2+\sum \limits_{k=0}^{N-1} \|{x(k)} \|_{Q}^2+\|{u(k)} \|_{R}^2 \label{eq:costFunction}\\
&\text{s.t.} 
  \ \ \, {x}(k+1)=A{x}(k)+B{u}(k), \ k=0, \ldots, N-1, \label{eq:Sys}\\
&\quad \ \ \ {x}(k) \in \mathcal{X}, \quad
{u}(k) \in \mathcal{U}, \quad \quad \ \ k=0, \ldots, N-1,\\
&\ \ \ \ \ \ \,  {x}(N) \in \mathcal{T}
\end{align}
\end{subequations}
that is periodically solved for a given initial condition $x(0)$ and horizon $N$, where $x(k) \in \R^n$ and $u(k) \in \R^m$ are the state and input variables at time step $k$, respectively. All matrices have the obvious dimensions and we assume stabilizability of the pair $(A,B)$, detectability of the pair $(Q^{\frac{1}{2}},A)$, $Q \succ 0$ and $R \succ 0$. Moreover, we assume $\mathcal{X}$, $\mathcal{U}$ and $\mathcal{T} \subseteq \mathcal{X}$ to be convex and compact polytopes that contain the origin as an interior point. 
The matrix $P$ is the solution to the discrete-time algebraic Riccati equation and $\mathcal{T}$ is the largest set on which the linear quadratic regulator satisfies the constraints and stabilizes the system. The terminal set $\mathcal{T}$  is determined with the procedure proposed in \cite{Gilbert1991}. 
Inserting \eqref{eq:Sys} into the cost function \eqref{eq:costFunction} yields the quadratic program (QP) 
\begin{align}
\begin{split}
&\min \limits_{\bar{U}} \ \frac{1}{2}\bar{U}'H\bar{U}+x(0)'F\bar{U}+\frac{1}{2}x(0)'Yx(0) \\ 
&\text{s.t.} \quad G\bar{U} \leq w+Ex(0),
\end{split}
\label{eq:ReformulatedMPCProblem}
\end{align} 
where $Y \in \R^{n \times n}$, $F \in \R^{n \times mN}$, $H \in \R^{mN \times mN}$, $G \in \R^{q \times mN}$, $w \in \R^q$, $E \in \R^{q \times n}$, $q$ is the number of constraints and $\bar{U}(x)=({u}(0)', \ldots, {u}(N-1)')'$. A state $x \in \R^n$ is \textit{feasible} if a control input $\bar{U}(x)$ exists that satisfies the constraints in \eqref{eq:ReformulatedMPCProblem}. Note that the positive definiteness of $H$, which is a consequence of $Q \succ 0$, $R \succ 0$ and $P \succ 0$, implies a unique optimal input sequence 
\begin{align}\label{eq:uniqueInputSequence}
\bar{U}^\star=(u^{\star}(0)', u^{\star}(1)', \ldots, u^{\star}(N-1)')' 
\end{align}
for every $x \in \mathcal{X}_f$, where $\mathcal{X}_f$ refers to the set of initial states for which problem \eqref{eq:ReformulatedMPCProblem} has a solution. 
We introduce the sets of active and inactive constraints
\begin{align}
\begin{split}
\mathcal{A}(x)&=\{i \in \mathcal{Q}~|~G^i\bar{U}^\star-w^i-E^ix=0 \},\\
\mathcal{I}(x)&=\{i \in \mathcal{Q}~|~G^i\bar{U}^\star-w^i-E^ix<0 \}= \mathcal{Q}\backslash\mathcal{A}(x),\\
\end{split}
\label{eq:Sets}
\end{align}
where $\mathcal{Q}:=\{1, \ldots, q \}$ denotes the set of all constraint indices. 
The following lemma, which is simple but central to all regional approaches~(\cite{Jost2015a, Koenig2017d, Koenig2017c, Koenig2017a, SchulzeDarup2017CDC}), is based on Theorem 2 in \cite{Bemporad2002}.                                                                                                                                                                                                                       
\begin{lemma}\label{lem:Jost}
Let $x \in \mathcal{X}_f$ be arbitrary and $\mathcal{A}(x)=\mathcal{A}$ the corresponding active set. Assume the matrix $G^{\mathcal{A}}$ has full row rank. Let 
\begin{align}\label{eq:LemmaJost}
  \begin{split}
    \bar{K}^\star&=H^{-1}(G^{\mathcal{A}})'(G^{\mathcal{A}}H^{-1}(G^{\mathcal{A}})')^{-1}S^{\mathcal{A}}-H^{-1}F',\\
    \bar{b}^\star&=H^{-1}(G^{\mathcal{A}})'(G^{\mathcal{A}}H^{-1}(G^{\mathcal{A}})')^{-1}w^{\mathcal{A}},\\
     T^\star&=\begin{pmatrix}G^{\mathcal{I}}H^{-1}(G^{\mathcal{A}})'(G^{\mathcal{A}}H^{-1}(G^{\mathcal{A}})')^{-1}S^{\mathcal{A}}-S^{\mathcal{I}} \\ (G^{\mathcal{A}}H^{-1}(G^{\mathcal{A}})')^{-1}S^{\mathcal{A}} \end{pmatrix},\\
      d^\star&=-\begin{pmatrix}G^{\mathcal{I}}H^{-1}(G^{\mathcal{A}})'(G^{\mathcal{A}}H^{-1}(G^{\mathcal{A}})')^{-1}w^{\mathcal{A}}-w^{\mathcal{I}} \\ (G^{\mathcal{A}}H^{-1}(G^{\mathcal{A}})')^{-1}w^{\mathcal{A}} \end{pmatrix}
  \end{split}
\end{align}
where  $S=E+GH^{-1}F'$, $S \in \R^{q \times n}$ and let $ \mathcal{P}^\star=\{ x \in \R^n ~|~T^\star x \leq d^\star \}$. Then, $\bar{U}^\star(\cdot): \mathcal{P}^\star \rightarrow \R^{Nm}$ defined by $x \rightarrow \bar{K}^\star x+\bar{b}^\star$
is the optimal input sequence $\bar{U}^\star \in \R^{mN}$ according to \eqref{eq:uniqueInputSequence} for all $x \in \mathcal{P}^\star$. Moreover, ${u}^\star (\cdot) : \mathcal{P}^\star \rightarrow \R^m$, i.e. the first $m$ components of $\bar{U}^\star (\cdot)$, defined by $x \rightarrow {K}^\star x+{b}^\star$ with $K^\star=\bar{K}^{\star \mathcal{M}}$ and $b^\star=\bar{b}^{\star \mathcal{M}}$ with $\mathcal{M}=\{1, \ldots m\}$ is the optimal feedback.
\end{lemma}
Lemma \ref{lem:Jost} yields an \textit{affine control law} and its polytope of validity from the solution of \eqref{eq:MPCProblem} for a \textit{single point} $x$. This result is used in a regional MPC approach that has been proposed in \cite{Jost2015a} and can be described as follows: Problem \eqref{eq:ReformulatedMPCProblem} is solved for the current state $x \in \mathcal{X}_f$. 
With the optimal input sequence $\bar{U}^\star \in \R^{mN}$ an active and inactive set according to \eqref{eq:Sets} results. This set is used to compute an optimal affine control law and its polytope containing the current state according to Lemma \ref{lem:Jost}. Instead of solving a QP in the next time step the feedback law, i.e., the first $m$ elements of the control law, are applied to the system. Only when leaving the current polytope, a new QP has to be solved to determine the next feedback law and polytope. 

\section{Regional MPC with suboptimally extended regions of validity}\label{sec:mainpart}
The optimal regional MPC approach avoids to solve a QP in every time step by reusing a feedback law on its polytope $\mathcal{P}^\star$ according to Lemma \ref{lem:Jost}. In this section we derive a region for the feedback law that is larger than $\mathcal{P}^\star$ by dropping optimality. 
Consequently, the number of QPs that must be solved are reduced further compared to the optimal approach, while states and inputs still satisfy the constraints. 
The new regions arise from intersecting a region required for feasibility and a region required for stability (see Fig. \ref{fig:idea} for a sketch), which are introduced in Sections \ref{sec:FeasibilityRegion} and \ref{sec:StabilityRegion} below. The intersection constitutes a nonlinearly bounded region, which is introduced in Section \ref{sec:NonlinearlyBoundedRegion} and used in an event-triggered approach in Section \ref{sec:NonlinearlyBoundedRegion}.    

\subsection{Feasibility region}\label{sec:FeasibilityRegion}

\begin{figure}[]
\footnotesize
\centering
{\begin{overpic}[scale=0.55,tics=10]%
{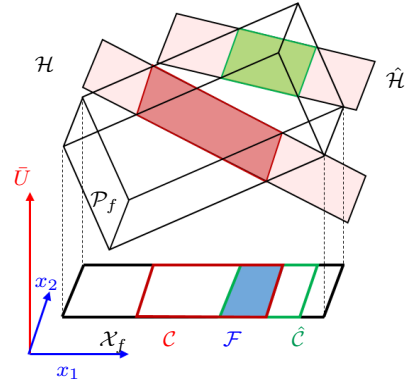}
\put(5,50){$\textcolor{red}{\bar{U}}$}
\put(10,77){$\mathcal{H}$}
\put(92,73){$\hat{\mathcal{H}}$}
\put(25,12){$\mathcal{X}_f$}
\put(23,45){$\mathcal{P}_f$}
\put(40,12){$\textcolor{red}{\mathcal{C}}$}
\put(54,12){$\textcolor{blue}{\mathcal{F}}$}
\put(70,12){$\textcolor{ForestGreen}{\hat{\mathcal{C}}}$}
\put(15,5){$\textcolor{blue}{x_1}$}
\put(10,26){$\textcolor{blue}{x_2}$}
\end{overpic}}   
\captionsetup{width=1\linewidth}
\caption{In this sketch, $\bar{K}^\star x + \bar{b}^\star$ according to Lemma \ref{lem:Jost} consists of two affine sub-laws defining the inputs $u(0)$ and $u(1)$, respectively. The constraints in \eqref{eq:ReformulatedMPCProblem} describe a polytope $\mathcal{P}_f$ in the augmented input space (white box). In this space, the affine laws associated with $u(0)$ and $u(1)$ define linear subspaces $\mathcal{H}$ and $\hat{\mathcal{H}}$, respectively (light red planes). 
The intersections $\mathcal{P}_f \cap \mathcal{H}$ and $\mathcal{P}_f \cap \hat{\mathcal{H}}$ are projected onto the $x$-space, which results in $\mathcal{C}$ and $\hat{\mathcal{C}}$. The intersection $\mathcal{C} \cap \hat{\mathcal{C}}$ yields a feasibility region $\mathcal{F}$ for the complete input sequence resulting from $\bar{K}^\star x + \bar{b}^\star$ (blue region). The feasibility region $\mathcal{C}$ for the input $u(0)$ resulting from $K^\star x + b^\star$ is larger than the feasibility region $\mathcal{F}$.} \label{fig:projections} 
\end{figure}

\begin{figure*}
\footnotesize
\quad 
\subfloat[\label{fig:polytope} Two feasibility regions for the feedback law that can be calculated online from the pointwise solution. The feasibility region for the control law according to \eqref{eq:FeasibilityN} is larger than the optimal polytope according to \eqref{eq:LemmaJost}. ]{\begin{overpic}[scale=0.55,tics=10]%
{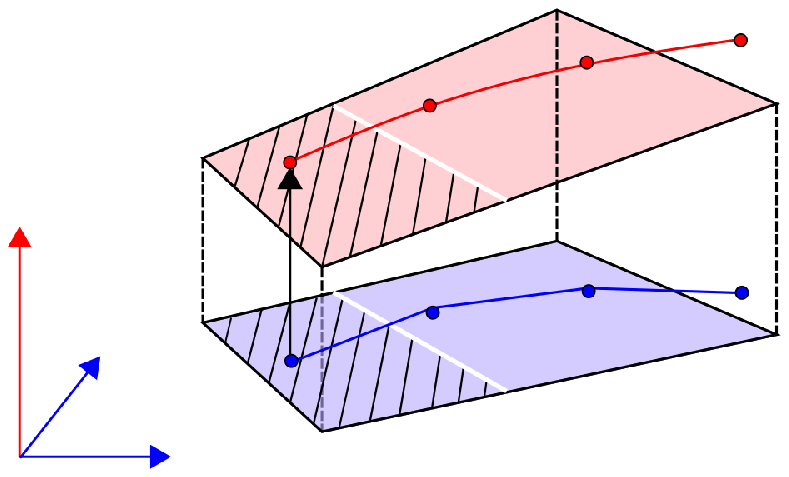}
\put(1,35){$\color{red} \bar{U}$}
\put(12,18){$ \color{blue} x_1$}
\put(14,37){$ \color{blue} x_2$}
\put(30,43){\footnotesize  QP}
\end{overpic}} \quad \quad  \quad 
\subfloat[\label{fig:quadric} The feedback law leads to a decreasing cost function, i.e., is stabilizing, on a region that can be described by a simple quadratic inequality.]{\begin{overpic}[scale=0.55,tics=10]%
{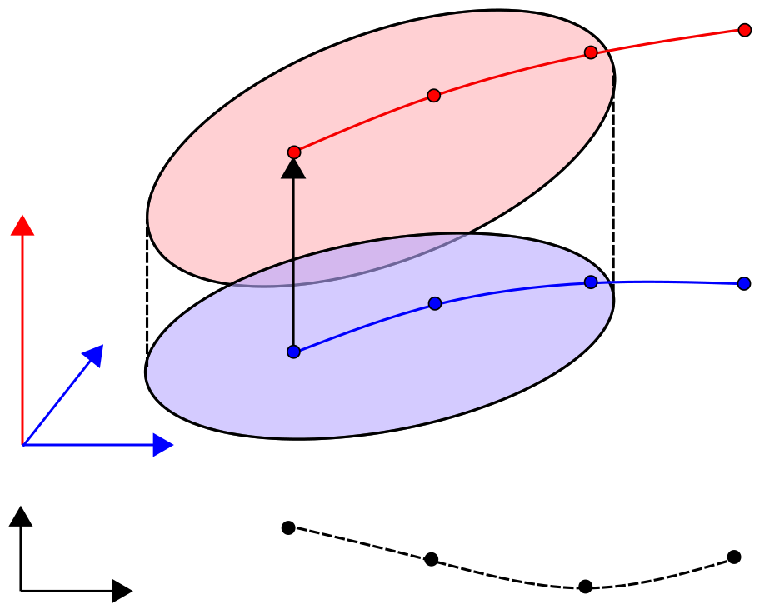}
\put(1,35){$ \color{red} \bar{U}$}
\put(12,19){$ \color{blue} x_1$}
\put(14,38){$ \color{blue} x_2$}
\put(0,9){$ V $}
\put(11,0){$ k $}
\put(29,34){\footnotesize  QP}
\end{overpic}} \quad \quad  \quad 
\subfloat[\label{fig:region} The feedback law can be reused as long as the system state remains in a nonlinearly bounded region. This region results from the intersection of the feasibility (see \ref{fig:polytope}) and stability region (see \ref{fig:quadric}).]{\begin{overpic}[scale=0.55,tics=10]%
{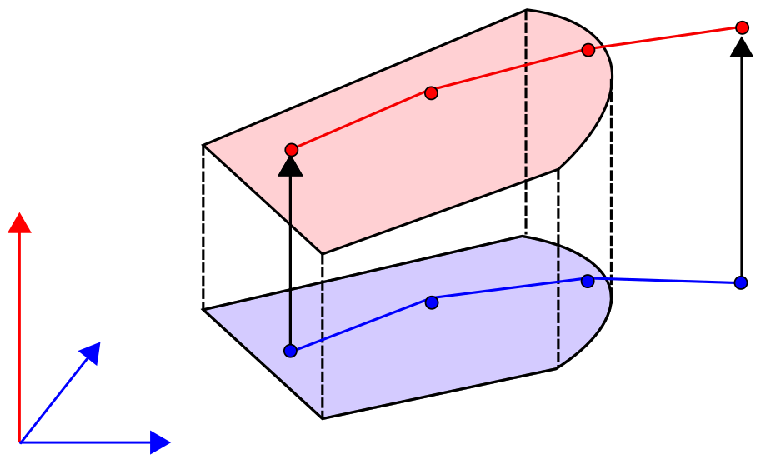}
\put(1,35){$ \color{red} \bar{U}$}
\put(12,18){$ \color{blue} x_1$}
\put(14,37){$ \color{blue} x_2$}
\put(30,43){\footnotesize  QP}
\put(83,52){\footnotesize  QP}
\end{overpic}}
\captionsetup{width=1\linewidth}
\caption{ A nonlinearly bounded region of validity resulting from the intersection of a feasibility and stability region \cite{Koenig2017d,Koenig2017c}.} \label{fig:idea} 
\end{figure*}

We derive two regions $\mathcal{C}$ and $\mathcal{F}$ on which the feedback signal $u(0)$ resulting from $K^\star x + b^\star$ satisfies the constraints and feasibility of the following state $x^+=Ax+bu(0)$ is guaranteed. Since we drop optimality both regions are larger than the optimal polytope $\mathcal{P}^\star$ according to Lemma \ref{lem:Jost}. They differ with respect to their size and their computational complexity. The regions are defined in Proposition \ref{prop:Cl} and Proposition \ref{prop:Feasibility} below and illustrated in Fig. \ref{fig:projections}. 

For preparation of Proposition \ref{prop:Cl}, we introduce the augmented vector 
\begin{align*}
\epsilon=(x(0)', u(0)', u(1)', \ldots u(N-1)')'. 
\end{align*}
Then the constraints in \eqref{eq:ReformulatedMPCProblem} can be represented as a polytope
\begin{align}\label{eq:PolytopeConstraints}
\mathcal{P}_\mathcal{F} = \{ \epsilon\in \R^{n+mN}  \ | \ \begin{pmatrix} -E & G \end{pmatrix} \epsilon \leq w \}
\end{align}
in the augmented input space $\R^n \times \R^{mN}$. 
The feedback law $u(0)=K^\star x + b^\star$ according to Lemma \ref{lem:Jost} can also be represented in the augmented input space as a linear subspace 
\begin{align}\label{eq:SubspaceControlLaw}
\mathcal{H}= \{ \epsilon\in \R^{n+mN}  \ | \ \begin{pmatrix}  -{K}^\star & L \end{pmatrix} \epsilon= b^{\star} \},
\end{align}
where $L \in \R^{m \times mN}$. 
\begin{proposition}\label{prop:Cl}
Let $x \in \mathcal{X}_f$ be an arbitrary state and $\mathcal{A}(x)$ the corresponding active set. Assume $G^\mathcal{A}$ has full row rank and let ${K}^\star$ and ${b}^\star$ be defined according to Lemma \ref{lem:Jost}. Let
\begin{align}\label{eq:FeasibilitySubLaw}
\mathcal{C}=\text{proj}_x(\mathcal{P}_{\mathcal{F}} \cap \mathcal{H}),
\end{align}
where $\mathcal{P}_{\mathcal{F}}$ and $\mathcal{H}$ are defined according to \eqref{eq:PolytopeConstraints} and \eqref{eq:SubspaceControlLaw}, respectively. Then for all $x \in \mathcal{C}$ the feedback law $K^\star x+b^\star$ yields a control input $u(0) \in \mathcal{U}$ that satisfies the constraints. Moreover, the following state $x^+=Ax+Bu(0)$ is feasible and $\mathcal{C} \supseteq \mathcal{P}^\star$ holds for the polytope $\mathcal{P}^\star$ from Lemma \ref{lem:Jost}.  
\end{proposition}

\begin{proof} By the definitions of $\mathcal{P}_f$ according to \eqref{eq:PolytopeConstraints} and $\mathcal{H}$ according to \eqref{eq:SubspaceControlLaw}, the intersection $\mathcal{P}_\mathcal{F} \cap \mathcal{H}$ yields all pairs $(x,\bar{U})$ with $\bar{U}=(u(0)',u(1)',\ldots,u(N-1)')'$ that fulfill the conditions $G\bar{U}-Ex \leq w$ and $u(0)=K^\star x+b^\star$. 
The states of these pairs can be determined by projecting the solution space onto the $x$-space, i.e. $\mathcal{C}=\text{proj}_x(\mathcal{P}_\mathcal{F} \cap \mathcal{H})$. Then, for all $x \in \mathcal{C}$, the feedback law $K^\star x+b^\star$ yields a control input $u(0) \in \mathcal{U}$, which proves the first part of Proposition \ref{prop:Cl}. For the second part, we use the fact that a feasible input sequence $\bar{U}$ for $x \in \mathcal{C}$ exists by construction. Then a feasible input sequence for $x^+=Ax+Bu(0)$ can be constructed by discarding the first $m$ elements of $\bar{U}$ and adding inputs resulting from the terminal controller. The relation $\mathcal{C} \supseteq \mathcal{P}^\star$ is obvious, because $\mathcal{C}$ guarantees constraint satisfaction only for $u(0)$ resulting from $K^\star x + b^\star$ and not the remaining inputs resulting from $\bar{K}^\star x + \bar{b}^\star$. Moreover, $\mathcal{C}$ does not take optimality into account.   
\hfill $\square$
\end{proof}

The region $\mathcal{C}$ in Proposition \ref{prop:Cl} is the largest region on which the input $u(0)$ resulting from $K^\star x + b^\star$ satisfies the constraints while feasibility of the following state $x^+$ is ensured. We emphasize that a large region is desirable because a feedback law may be reused for more time steps than in the existing approach, thus reducing the number of QPs to be solved (see Fig. \ref{fig:idea}). However, a computational expensive projection is required. Consequently, the regions $\mathcal{C}$ have to be computed offline before. They can be used in our regional MPC approach at runtime. The region  $\mathcal{F}$ introduced in Proposition \ref{prop:Feasibility} below is smaller than the region $\mathcal{C}$ from Proposition \ref{prop:Cl}, but no projections are needed so that these regions can be computed at runtime. The region $\mathcal{F}$ ensures that the input sequence $\bar{U}=(u(0)',u(1)', \ldots, u(N-1)')'$ resulting from $\bar{K}^\star x + \bar{b}^\star$ satisfies the constraints. Note that this is different from $\mathcal{C}$ because all inputs $u(0),u(1), \ldots, u(N-1)$ have to be defined by the law $\bar{K}^\star x + \bar{b}^\star$. The region $\mathcal{F}$ can be calculated without any projection by inserting the affine law $\bar{K}^\star x + \bar{b}^\star$ into the constraints in \eqref{eq:ReformulatedMPCProblem}. The region $\mathcal{F}$ is characterized formally in the following proposition that is based on the results in \cite{Bemporad2001}. 
\begin{proposition}\label{prop:Feasibility}
Let $x \in \mathcal{X}_f$ be arbitrary and $\mathcal{A}(x)=\mathcal{A}$ the corresponding active set. Assume the matrix $G^{\mathcal{A}}$ has full row rank and let $\bar{K}^\star$ and $\bar{b}^\star$ be defined as in Lemma \ref{lem:Jost}. Let
\begin{align}\label{eq:FeasibilityN}
\begin{split}
T_{1}&=G^{\mathcal{I}}H^{-1}(G^{\mathcal{A}})'(G^{\mathcal{A}}H^{-1}(G^{\mathcal{A}})')^{-1}S^{\mathcal{A}}-S^{\mathcal{I}}, \\
d_{1}&=-G^{\mathcal{I}}H^{-1}(G^{\mathcal{A}})'(G^{\mathcal{A}}H^{-1}(G^{\mathcal{A}})')^{-1}w^{\mathcal{A}}+w^{\mathcal{I}},
\end{split}
\end{align}
where  $S=E+GH^{-1}F'$, $S \in \R^{q \times n}$ and let $ \mathcal{F}=\{ x \in \R^n ~|~T_1x \leq d_1 \}$. Then, for all $x \in \mathcal{F}$, the affine law $\bar{K}^\star x+\bar{b}^\star$ yields a control input sequence $\bar{U}$ that satisfies the constraints in \eqref{eq:ReformulatedMPCProblem}. Moreover, the successor state $x^+=Ax+Bu(0)$ is feasible and $\mathcal{F}\supseteq \mathcal{P}^\star$ holds for the polytope $\mathcal{P}^\star$ from Lemma \ref{lem:Jost}.
\end{proposition}
A proof for Proposition \ref{prop:Feasibility} can be found in \cite{Koenig2017c}.

\subsection{Stability region}\label{sec:StabilityRegion}
In order to ensure stability, the cost function has to decrease along all closed-loop trajectories. More precisely, we need to enforce
\begin{align}\label{eq:CondStability}
V(x(k),\bar{U}(x(k))) < \lambda V(x(k-1),\bar{U}(x(k-1)))
\end{align}
with an arbitrary but fixed $\lambda \in (0,1]$. Note that on the feasibility regions $\mathcal{C}$ and $\mathcal{F}$ according to Proposition \ref{prop:Cl} and \ref{prop:Feasibility} the feedback law can lead to an increasing cost function. To avoid this, we have to calculate a second region on which stability is guaranteed and intersect it with the feasibility region. 
A stability region can be obtained by substituting the law $\bar{K}^\star x + \bar{b}^\star$ and the system dynamics into \eqref{eq:CondStability}. As a result we obtain a quadratic inequality that depends only on the current state $x$ and describes a nonlinearly bounded state-space region $\mathcal{V}$. This is stated more precisely in Proposition~\ref{prop:Stability}. In preparation of Proposition~\ref{prop:Stability} we introduce 
  \begin{equation}\label{eq:HilfsgroessenFuerV}
    \begin{split}
      T_{2}&=M_2-\lambda(2 M_4' M_1 M_3-M_2 M_3), \\
      T_{3}&=M_1-\lambda M_3'M_1M_3, \\
      d_{2}&=\lambda(M_4' M_1 M_4 + M_2 M_4+M_5)-M_5
      \end{split}
  \end{equation}
with $M_1=\frac{1}{2} \bar{K}^{\star \prime} H \bar{K}^\star + F \bar{K}^\star + \frac{1}{2} Y$, $M_2=\bar{b}^{\star \prime} H \bar{K}^\star + \bar{b}^{\star \prime}  F'$, $ M_3=(A+BK^{\star})^{-1}$ and $M_4=- M_3 B b^{\star}$.
\begin{proposition}\label{prop:Stability}
Let $\bar{K}^\star x + \bar{b}^\star$ be an arbitrary law according to Lemma \ref{lem:Jost}. Let $x^- \in \mathcal{X}_f$ be arbitrary and $x=Ax^- + Bu(0)$ be the successor state that results from applying this law to system \eqref{eq:Sys}. Let $\mathcal{V}$ be defined by
\begin{align}\label{eq:RegionStability}
\mathcal{V}= \left\{\zeta\in\R^n | \zeta^\prime T_3 \zeta+ T_2\zeta< d_2\right\}.
\end{align}
If $(A+B K^{\star}_0)^{-1}$ exists, then 
\begin{align}
V(x, \bar{U}(x)) < \lambda V(x^-, \bar{U}(x^-)) \ \Longleftrightarrow \ x \in \mathcal{V}
\end{align}
holds.
\end{proposition}
The proof of Proposition \ref{prop:Stability} can be found in \cite{Koenig2017c}. Note that the previous state $x^-$ is not required to be located within the optimal polytope $\mathcal{P}^\star$ so that the reusability of the feedback law is guaranteed not just in the first time step after leaving the current polytope but for later time steps, too. 
If the inverse and thus the region $\mathcal{V}$ cannot be determined, the optimal region $\mathcal{P}^\star$ according to \eqref{eq:LemmaJost} must be computed instead. Moreover, note that the decrease of the cost function and thus the size of the region $\mathcal{V}$ can be tuned with the parameter $\lambda$. By means of that the closed-loop performance can be varied.      

\subsection{Nonlinearly bounded regions of validity}\label{sec:NonlinearlyBoundedRegion}
Let us consider the intersection
\begin{align}\label{eq:RegionIntersection}
\mathcal{E}=\mathcal{B} \cap \mathcal{V},
\end{align}   
where $\mathcal{B}=\mathcal{C}$ according to Proposition \ref{prop:Cl} or $\mathcal{B}=\mathcal{F}$ according to Proposition \ref{prop:Feasibility} and $\mathcal{V}$ is defined according to Proposition \ref{prop:Stability}. Then obviously for all $x \in \mathcal{E}$ the feedback law  $K^\star x + b^\star$ yields an input $u(0)$ that satisfies the constraints and leads to a decreasing cost function. The region $\mathcal{E}$ generally is a nonlinearly bounded region since it is the intersection of a polytope and a region described by a quadratic inequality. Consequently, $\mathcal{E}$ may be a non-convex region, e.g., an $n$-dimensional hyperbola.
Nevertheless checking if a state $x \in \mathcal{X}_f$ is located inside the region $\mathcal{E}$ requires only a low effort since it merely requires to evaluate the polytope $T_1 x \leq d_1$ and the simple inequality $x'T_3x+T_2x <d_2$. We emphasize that we do not calculate the borders of the region $\mathcal{E}$  explicitly.   

\begin{figure}[]
\centering
\begin{overpic}[width=0.45\textwidth,tics=10]%
{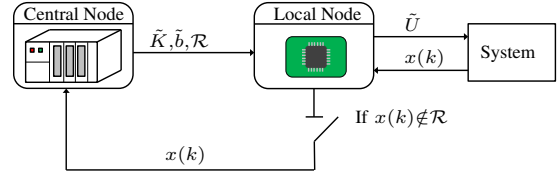}
\put(48,28){\scriptsize Local Node}
\put(86,21){\scriptsize System}
\put(2,28){\scriptsize Central Node}
\put(72,26){$\scriptstyle \tilde{U}$}
\put(25,23){$\scriptstyle \tilde{K}, \tilde{b}, \mathcal{R}$}
\put(72,20){$\scriptstyle x(k)$}
\put(63,9){$\scriptstyle \text{If} \ x(k) \notin \mathcal{R}$}
\put(28,2){$\scriptstyle x(k)$}
\end{overpic}
\caption{Event-triggered update of the stabilizing feedback law.}
\label{fig:NetworkedSetting}
\end{figure}  

\subsection{Event-triggered feedback law update}\label{sec:eventTriggered}
The feedback law $K^\star x+ b^\star$ and the new nonlinearly bounded regions $\mathcal{E}$ can be used in a regional MPC algorithm summarized as follows: For an arbitrary state $x \in \mathcal{X}_f$ the corresponding active set $\mathcal{A}(x)$ can be determined by solving the QP \eqref{eq:ReformulatedMPCProblem}. The active set can be used to calculate the feedback law $K^\star (x)+b^\star$ with Lemma \ref{lem:Jost} and its region of validity $\mathcal{R}$. The region $\mathcal{R}=\mathcal{E}$ is computed with \eqref{eq:RegionIntersection} if the inverse $(A+BK^{\star})^{-1}$ exists. 
Otherwise, the optimal polytope $\mathcal{R}= \mathcal{P}^\star$ is calculated according to Lemma \ref{lem:Jost}. As long as the system remains in the region $\mathcal{R}$, the feedback law can be reused. Leaving the region $\mathcal{R}$ is the event that triggers solving a new QP to determine the next active set and thus the next feedback law and validity region. The procedure described above leads to the time varying feedback law 
\begin{align}\label{eq:TimeVaryingControlLaw}
\tilde{U}(x)&=\tilde{K}(k)x(k)+\tilde{b}(k) 
\end{align}
and the region $\mathcal{R}(k)$ with $\mathcal{R}(0)=\mathcal{R}_{new}(0)$, $\tilde{K}(0)={K}^\star$, $\tilde{b}(0)=b^\star$ and
\begin{align}\label{eq:NewUpdateRule2}
\begin{pmatrix}
\tilde{K}(k) \\
\tilde{b}(k) \\
\mathcal{R}(k)
\end{pmatrix}=\begin{cases}
\begin{pmatrix}
\tilde{K}(k-1) \\
\tilde{b}(k-1) \\
\mathcal{R}(k-1)
\end{pmatrix} & \begin{matrix}  \text{if} \ x(k) \in \mathcal{R} (k-1) \end{matrix}\\
\begin{pmatrix}
{K}^\star \\
{b}^\star\\
\mathcal{R}_{new}
\end{pmatrix} &\text{otherwise}
\end{cases}
\end{align}
for all $k>0$, where
\begin{align}\label{eq:RegionDecision}
\mathcal{R}_{new}(k)=\begin{cases} 
\mathcal{E} \quad \text{if} \ (A+B \tilde{K}(k))^{-1} \ {\rm exists}   \\
\mathcal{P}^\star \quad \text{otherwise.}
\end{cases}
\end{align}

Applying the update rule \eqref{eq:NewUpdateRule2} to system \eqref{eq:Sys} results in an asymptotically stable closed-loop system with constraint satisfaction. This is stated more precisely in the following proposition. 

\begin{proposition}\label{ko:ZulaessigkeitStabilitaet}
Let $P$ and $\mathcal{T}$ be defined as in Sect. \ref{sec:problemstatement} and $\mathcal{X}_f$ the set of states, such that~\eqref{eq:ReformulatedMPCProblem} has a feasible solution. Let the control law \eqref{eq:TimeVaryingControlLaw}--\eqref{eq:RegionDecision} be applied to the system \eqref{eq:Sys}. Then the origin is an asymptotically stable steady state of the controlled system  with domain of attraction $\mathcal{X}_f$. 
\end{proposition}

The proof of Proposition \ref{ko:ZulaessigkeitStabilitaet} is stated in \cite{Koenig2017c}.

\subsection{Networked control setting}
The new event-triggered MPC approach is appropriate for a networked control setting (see Fig. \ref{fig:NetworkedSetting}). The setting is the same as for the optimal approach presented in \cite{BernerP2016}. Our new nonlinearly bounded regions \eqref{eq:RegionIntersection} can be evaluated with a low effort on a local node resulting in a lean local node as in the optimal approach. The procedure can be described as follows.

The matrices ($K^{\star}$, $b^{\star}$) defining the feedback law and the matrices ($T^\star$, $d^\star$) or ($T_1$, $T_2$, $T_3$, $d_1$, $d_2$) defining the corresponding region of validity $\mathcal{R}$ are calculated on a powerful central node by solving the QP \eqref{eq:ReformulatedMPCProblem}. The matrices are transmitted to a local node, where the feedback law is evaluated with a low computational effort. 
Closed-loop control signals are applied to the system as long as the current state remains in the corresponding region. If the current region is left and thus the feedback law becomes invalid, the local node requests a new feedback law and region from the central node. Compared to the optimal approach the effort on the central node can be reduced since the number of solved QPs can be reduced. 
The effort on the local node is similar to the effort in the optimal approach. It can be expressed by floating point operations (flops) per time step. Checking whether a point in state space lies in the region $\mathcal{R}$ in our new approach requires 
\begin{align}\label{eq:flopsNew}
\begin{split}
2qn \quad &\text{if}\quad \mathcal{R}=\mathcal{P}^\star, \\
  2n^2+3n+2\tilde{q}n \quad &\text{if} \quad \mathcal{R}=\mathcal{E},
  \end{split}
\end{align}    
where $\tilde{q}$ denotes the number of inactive constraints. The leading term in \eqref{eq:flopsNew} results from the evaluation of $\zeta' T_3 \zeta$ in \eqref{eq:RegionStability} that includes  $n$ inner products of length $n$. Since the region $\mathcal{R}$ does not depend on the control input sequence it can be evaluated with a low effort even for systems with a large number of optimization variables $mN$. 

\section{Examples}\label{sec:example}
In this section we apply the new approach to two examples that differ with respect to their sizes. For both systems we calculate closed-loop system trajectories for random initial states until $||x_k|| \leq 10^{-2}$. 

\begin{example}
We consider the system
\begin{align}
Y(s)=\frac{2}{s^2+s+2} U(s) \nonumber
\end{align}
that is discretized with a sample time $T_s=\unit[0.1]{s}$. The matrices of the system \eqref{eq:Sys} are
\begin{align}
A=\begin{pmatrix}  
0.8955 & -0.1897 \\
0.0948 & 0.9903
\end{pmatrix}, \quad
B=\begin{pmatrix}
0.0948 \\
0.0048
\end{pmatrix}. \nonumber
\end{align}
The example is similar to the one used by \cite{Seron2003}, but the states and inputs are subject to the constraints $-3 \leq x_i \leq 3, \ i=1,2$ and $-2 \leq u_i \leq 2, \ i=1$, respectively. The weighting matrices read $Q=\text{diag}(0.01, 4)$ and $R=0.01$. The weighting matrix $P$ and the terminal set $\mathcal{T}$ are determined according to Sect. \ref{sec:problemstatement}. With a prediction horizon of $N= 4$, the system results in a QP with 32 inequalities and 4 optimization variables.
\end{example}

\begin{example}
We consider the system 
\begin{align}
G(s)=\begin{pmatrix}
\frac{0.05}{36s^2+6s+1} & \frac{0.02(2s+1)}{8s+1}\\
\frac{0.02(2s+1)}{8s+1} & \frac{0.05}{12s^2+3s+1}
\end{pmatrix} \nonumber
\end{align}
that is discretized with a sample time $T_s=\unit[1]{s}$. A state-space model with $n=6$ and $m=2$ results. We refer to \cite{Jost2015b} for the state-space system matrices.
The state and input constraints read $-15 \leq x_i \leq 15$ for $i=1,\ldots,6$ and $-3 \leq u_j \leq 3$ for $j=1,2$. The states and inputs are weighted with $Q=10 I^{6 \times 6}$  and $R=0.01 I^{2 \times 2}$. The terminal state weighting matrix $P$ and the terminal set $\mathcal{T}$ are determined according to Sect. \ref{sec:problemstatement}. With a horizon of $N=40$, this example results in a QP with 670 inequalities and 80 optimization variables.
\end{example}

\begin{figure*}
\centering
\captionsetup[subfloat]{justification=centering}
\captionsetup[subfigure]{labelformat=empty}
\subfloat[]{
%
%
\begin{psfrags}%
\psfragscanon%
\scriptsize%
%
\psfrag{s02}[b][b]{\color[rgb]{0.15,0.15,0.15}\setlength{\tabcolsep}{0pt}\begin{tabular}{c}$u(k)$\end{tabular}}%
\psfrag{s04}[t][t]{\color[rgb]{0.15,0.15,0.15}\setlength{\tabcolsep}{0pt}\begin{tabular}{c}$k$\end{tabular}}%
\psfrag{s06}[b][b]{\color[rgb]{0.15,0.15,0.15}\setlength{\tabcolsep}{0pt}\begin{tabular}{c}$e(k)$\end{tabular}}%
\psfrag{s12}[b][b]{\color[rgb]{0.15,0.15,0.15}\setlength{\tabcolsep}{0pt}\begin{tabular}{c}$x(k)$\end{tabular}}%
%
\color[rgb]{0.15,0.15,0.15}%
%
\psfrag{x01}[t][t]{0}%
\psfrag{x02}[t][t]{2}%
\psfrag{x03}[t][t]{4}%
\psfrag{x04}[t][t]{6}%
\psfrag{x05}[t][t]{8}%
\psfrag{x06}[t][t]{10}%
\psfrag{x07}[t][t]{12}%
\psfrag{x08}[t][t]{14}%
\psfrag{x09}[t][t]{16}%
\psfrag{x10}[t][t]{18}%
\psfrag{x11}[t][t]{0}%
\psfrag{x12}[t][t]{2}%
\psfrag{x13}[t][t]{4}%
\psfrag{x14}[t][t]{6}%
\psfrag{x15}[t][t]{8}%
\psfrag{x16}[t][t]{10}%
\psfrag{x17}[t][t]{12}%
\psfrag{x18}[t][t]{14}%
\psfrag{x19}[t][t]{16}%
\psfrag{x20}[t][t]{18}%
\psfrag{x21}[t][t]{0}%
\psfrag{x22}[t][t]{2}%
\psfrag{x23}[t][t]{4}%
\psfrag{x24}[t][t]{6}%
\psfrag{x25}[t][t]{8}%
\psfrag{x26}[t][t]{10}%
\psfrag{x27}[t][t]{12}%
\psfrag{x28}[t][t]{14}%
\psfrag{x29}[t][t]{16}%
\psfrag{x30}[t][t]{18}%
%
\psfrag{v01}[r][r]{0}%
\psfrag{v02}[r][r]{}%
\psfrag{v03}[r][r]{1}%
\psfrag{v04}[r][r]{-2}%
\psfrag{v05}[r][r]{0}%
\psfrag{v06}[r][r]{2}%
\psfrag{v07}[r][r]{-1}%
\psfrag{v08}[r][r]{0}%
\psfrag{v09}[r][r]{1}%
%
\includegraphics[width = 0.28\textwidth]{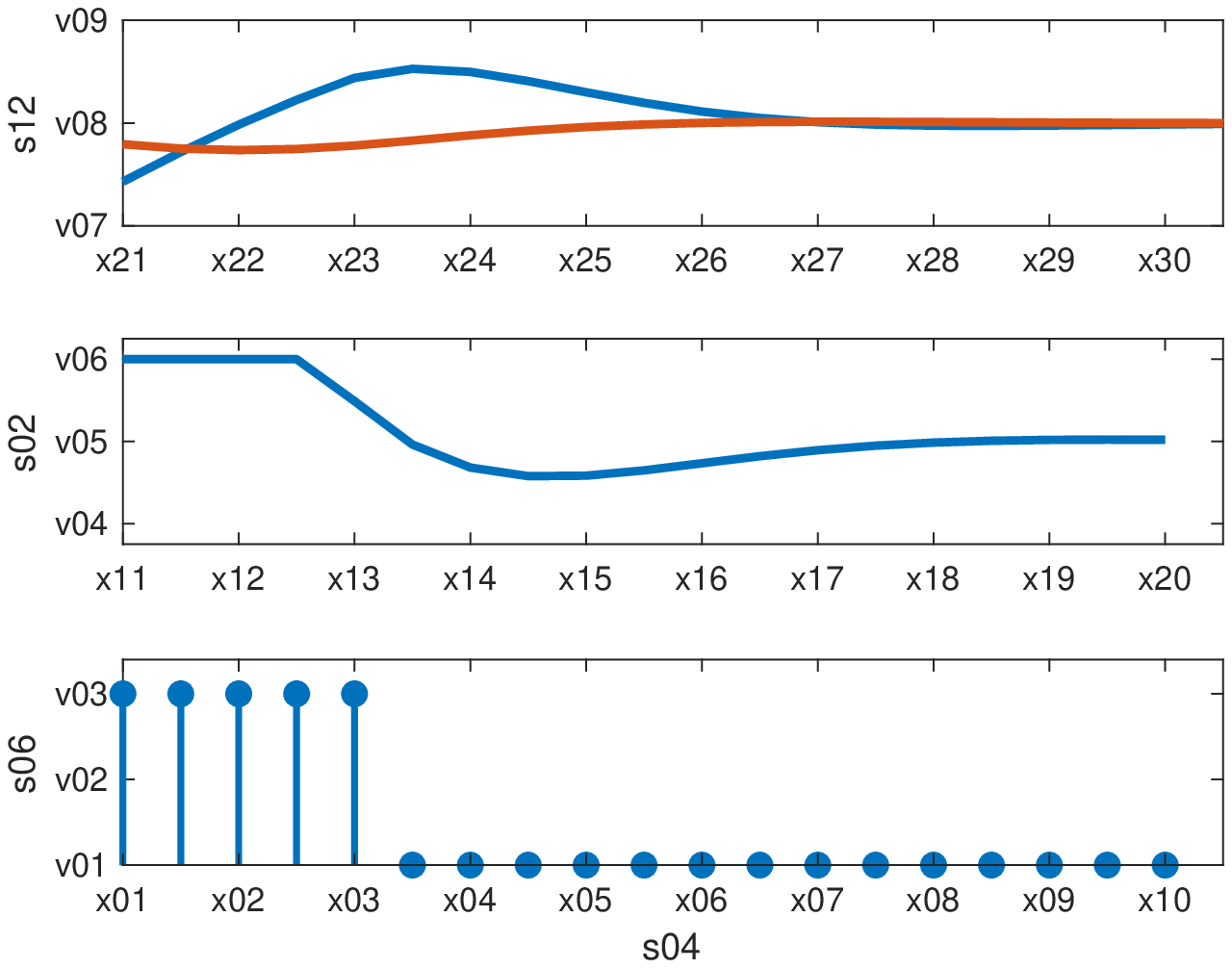}%
\end{psfrags}%
%
 \quad \quad 
%
%
\begin{psfrags}%
\psfragscanon%
\scriptsize%
%
\psfrag{s01}[b][b]{\color[rgb]{0.15,0.15,0.15}\setlength{\tabcolsep}{0pt}\begin{tabular}{c}$e(k)$\end{tabular}}%
\psfrag{s02}[b][b]{\color[rgb]{0.15,0.15,0.15}\setlength{\tabcolsep}{0pt}\begin{tabular}{c}$u(k)$\end{tabular}}%
\psfrag{s05}[t][t]{\color[rgb]{0.15,0.15,0.15}\setlength{\tabcolsep}{0pt}\begin{tabular}{c}$k$\end{tabular}}%
\psfrag{s06}[b][b]{\color[rgb]{0.15,0.15,0.15}\setlength{\tabcolsep}{0pt}\begin{tabular}{c}$x(k)$\end{tabular}}%
%
\color[rgb]{0.15,0.15,0.15}%
%
\psfrag{x01}[t][t]{0}%
\psfrag{x02}[t][t]{2}%
\psfrag{x03}[t][t]{4}%
\psfrag{x04}[t][t]{6}%
\psfrag{x05}[t][t]{8}%
\psfrag{x06}[t][t]{10}%
\psfrag{x07}[t][t]{12}%
\psfrag{x08}[t][t]{14}%
\psfrag{x09}[t][t]{16}%
\psfrag{x10}[t][t]{18}%
\psfrag{x11}[t][t]{0}%
\psfrag{x12}[t][t]{2}%
\psfrag{x13}[t][t]{4}%
\psfrag{x14}[t][t]{6}%
\psfrag{x15}[t][t]{8}%
\psfrag{x16}[t][t]{10}%
\psfrag{x17}[t][t]{12}%
\psfrag{x18}[t][t]{14}%
\psfrag{x19}[t][t]{16}%
\psfrag{x20}[t][t]{18}%
\psfrag{x21}[t][t]{0}%
\psfrag{x22}[t][t]{2}%
\psfrag{x23}[t][t]{4}%
\psfrag{x24}[t][t]{6}%
\psfrag{x25}[t][t]{8}%
\psfrag{x26}[t][t]{10}%
\psfrag{x27}[t][t]{12}%
\psfrag{x28}[t][t]{14}%
\psfrag{x29}[t][t]{16}%
\psfrag{x30}[t][t]{18}%
%
\psfrag{v01}[r][r]{0}%
\psfrag{v02}[r][r]{}%
\psfrag{v03}[r][r]{1}%
\psfrag{v04}[r][r]{-2}%
\psfrag{v05}[r][r]{0}%
\psfrag{v06}[r][r]{2}%
\psfrag{v07}[r][r]{-1}%
\psfrag{v08}[r][r]{0}%
\psfrag{v09}[r][r]{1}%
%
\includegraphics[width = 0.28\textwidth]{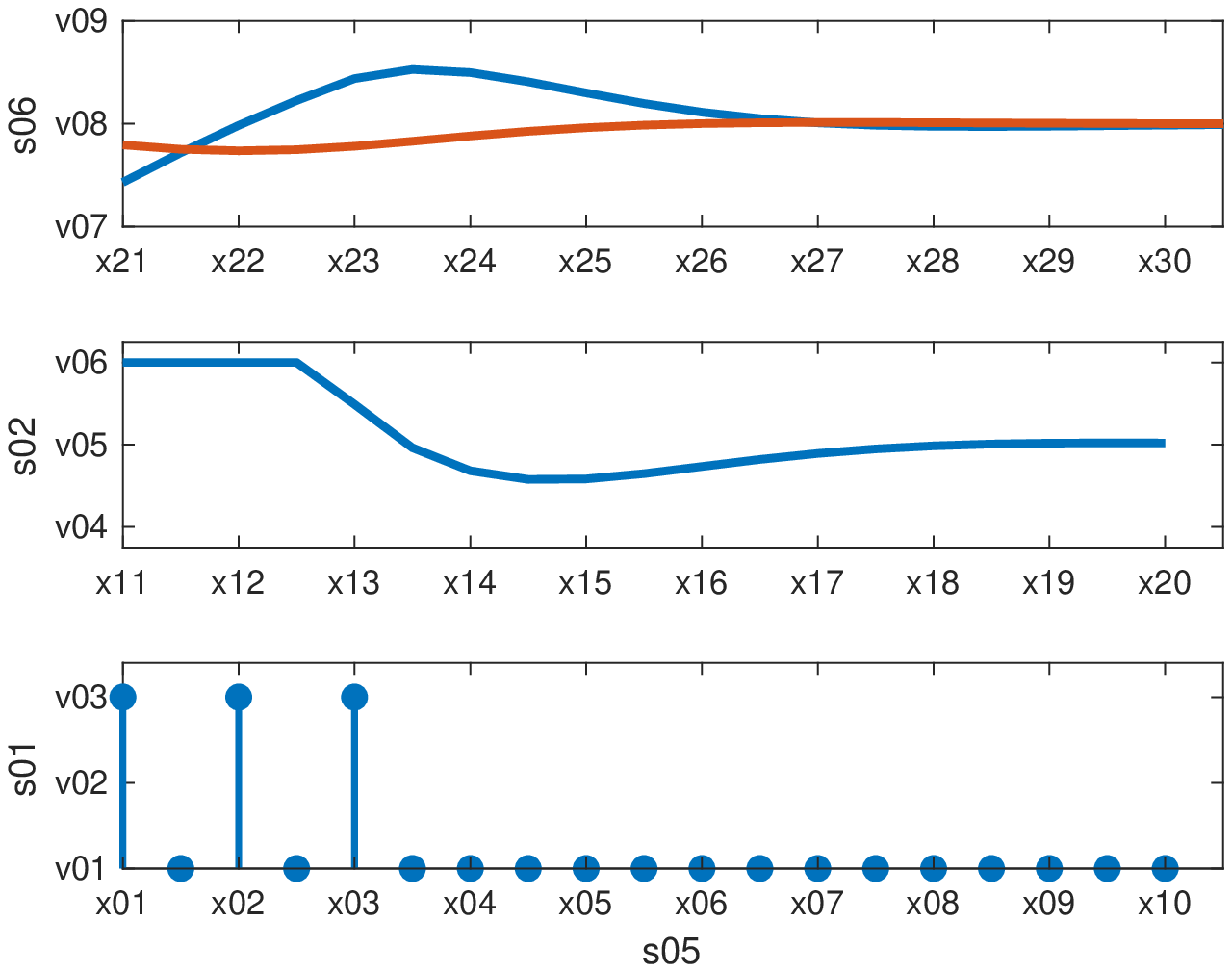}%
\end{psfrags}%
%
 \quad \quad 
%
%
\begin{psfrags}%
\psfragscanon%
\scriptsize%
%
\psfrag{s02}[b][b]{\color[rgb]{0.15,0.15,0.15}\setlength{\tabcolsep}{0pt}\begin{tabular}{c}$u(k)$\end{tabular}}%
\psfrag{s03}[t][t]{\color[rgb]{0.15,0.15,0.15}\setlength{\tabcolsep}{0pt}\begin{tabular}{c}$k$\end{tabular}}%
\psfrag{s07}[b][b]{\color[rgb]{0.15,0.15,0.15}\setlength{\tabcolsep}{0pt}\begin{tabular}{c}$x(k)$\end{tabular}}%
\psfrag{s12}[b][b]{\color[rgb]{0.15,0.15,0.15}\setlength{\tabcolsep}{0pt}\begin{tabular}{c}$e(k)$\end{tabular}}%
%
\color[rgb]{0.15,0.15,0.15}%
%
\psfrag{x01}[t][t]{0}%
\psfrag{x02}[t][t]{2}%
\psfrag{x03}[t][t]{4}%
\psfrag{x04}[t][t]{6}%
\psfrag{x05}[t][t]{8}%
\psfrag{x06}[t][t]{10}%
\psfrag{x07}[t][t]{12}%
\psfrag{x08}[t][t]{14}%
\psfrag{x09}[t][t]{16}%
\psfrag{x10}[t][t]{18}%
\psfrag{x11}[t][t]{0}%
\psfrag{x12}[t][t]{2}%
\psfrag{x13}[t][t]{4}%
\psfrag{x14}[t][t]{6}%
\psfrag{x15}[t][t]{8}%
\psfrag{x16}[t][t]{10}%
\psfrag{x17}[t][t]{12}%
\psfrag{x18}[t][t]{14}%
\psfrag{x19}[t][t]{16}%
\psfrag{x20}[t][t]{18}%
\psfrag{x21}[t][t]{0}%
\psfrag{x22}[t][t]{2}%
\psfrag{x23}[t][t]{4}%
\psfrag{x24}[t][t]{6}%
\psfrag{x25}[t][t]{8}%
\psfrag{x26}[t][t]{10}%
\psfrag{x27}[t][t]{12}%
\psfrag{x28}[t][t]{14}%
\psfrag{x29}[t][t]{16}%
\psfrag{x30}[t][t]{18}%
%
\psfrag{v01}[r][r]{0}%
\psfrag{v02}[r][r]{}%
\psfrag{v03}[r][r]{1}%
\psfrag{v04}[r][r]{-2}%
\psfrag{v05}[r][r]{0}%
\psfrag{v06}[r][r]{2}%
\psfrag{v07}[r][r]{-1}%
\psfrag{v08}[r][r]{0}%
\psfrag{v09}[r][r]{1}%
%
\includegraphics[width = 0.28\textwidth]{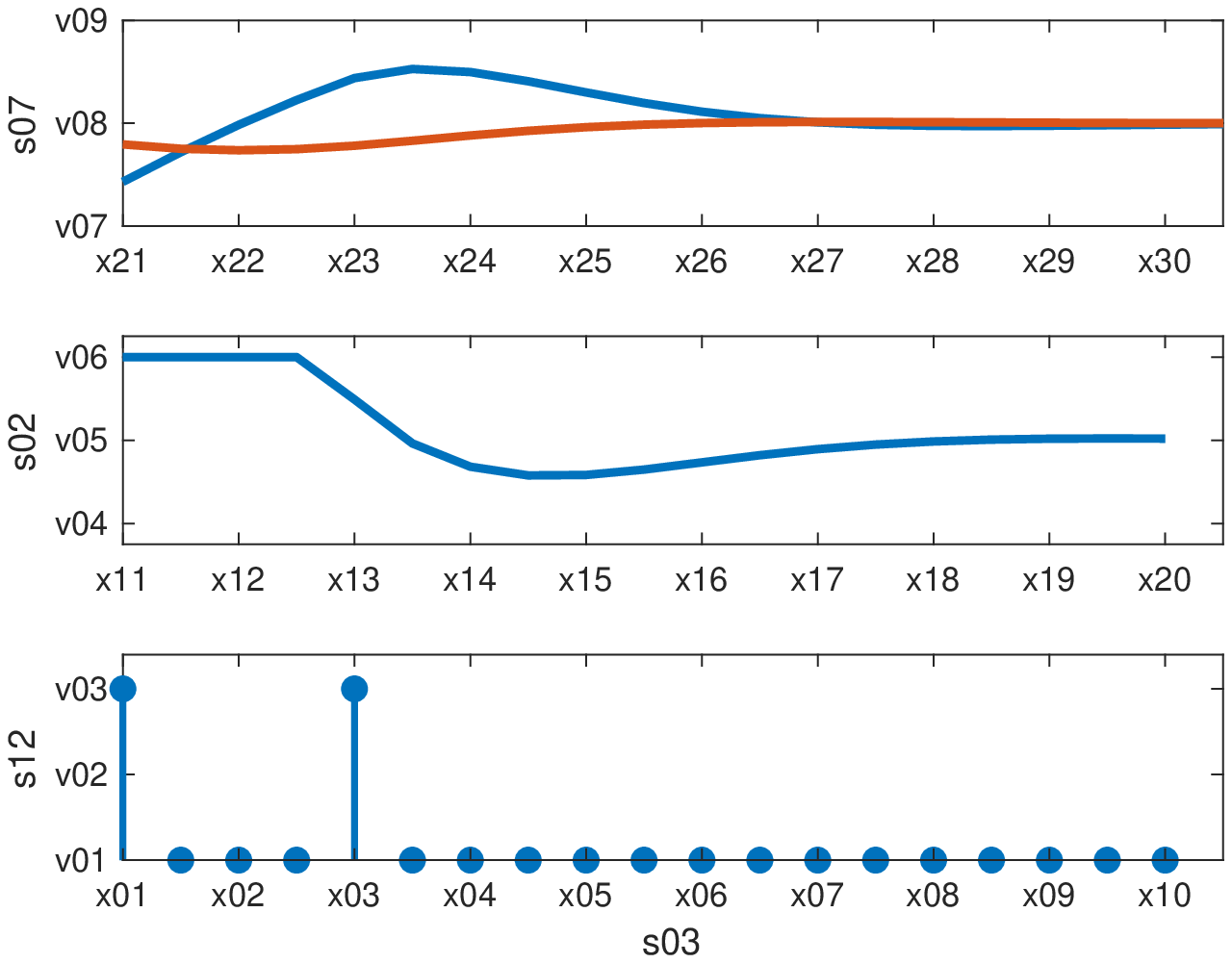}%
\end{psfrags}%
%
} \\ 
\subfloat[Optimal approach]{
%
%
\begin{psfrags}%
\psfragscanon%
\scriptsize%
%
\psfrag{s01}[b][b]{\color[rgb]{0.15,0.15,0.15}\setlength{\tabcolsep}{0pt}\begin{tabular}{c}$x_2$\end{tabular}}%
\psfrag{s04}[t][t]{\color[rgb]{0.15,0.15,0.15}\setlength{\tabcolsep}{0pt}\begin{tabular}{c}$x_1$\end{tabular}}%
%
\color[rgb]{0.15,0.15,0.15}%
%
\psfrag{x01}[t][t]{-2}%
\psfrag{x02}[t][t]{-1.5}%
\psfrag{x03}[t][t]{-1}%
\psfrag{x04}[t][t]{-0.5}%
\psfrag{x05}[t][t]{0}%
\psfrag{x06}[t][t]{0.5}%
\psfrag{x07}[t][t]{1}%
\psfrag{x08}[t][t]{1.5}%
\psfrag{x09}[t][t]{2}%
%
\psfrag{v01}[r][r]{-1}%
\psfrag{v02}[r][r]{-0.5}%
\psfrag{v03}[r][r]{0}%
\psfrag{v04}[r][r]{0.5}%
\psfrag{v05}[r][r]{1}%
%
\includegraphics[width = 0.28\textwidth]{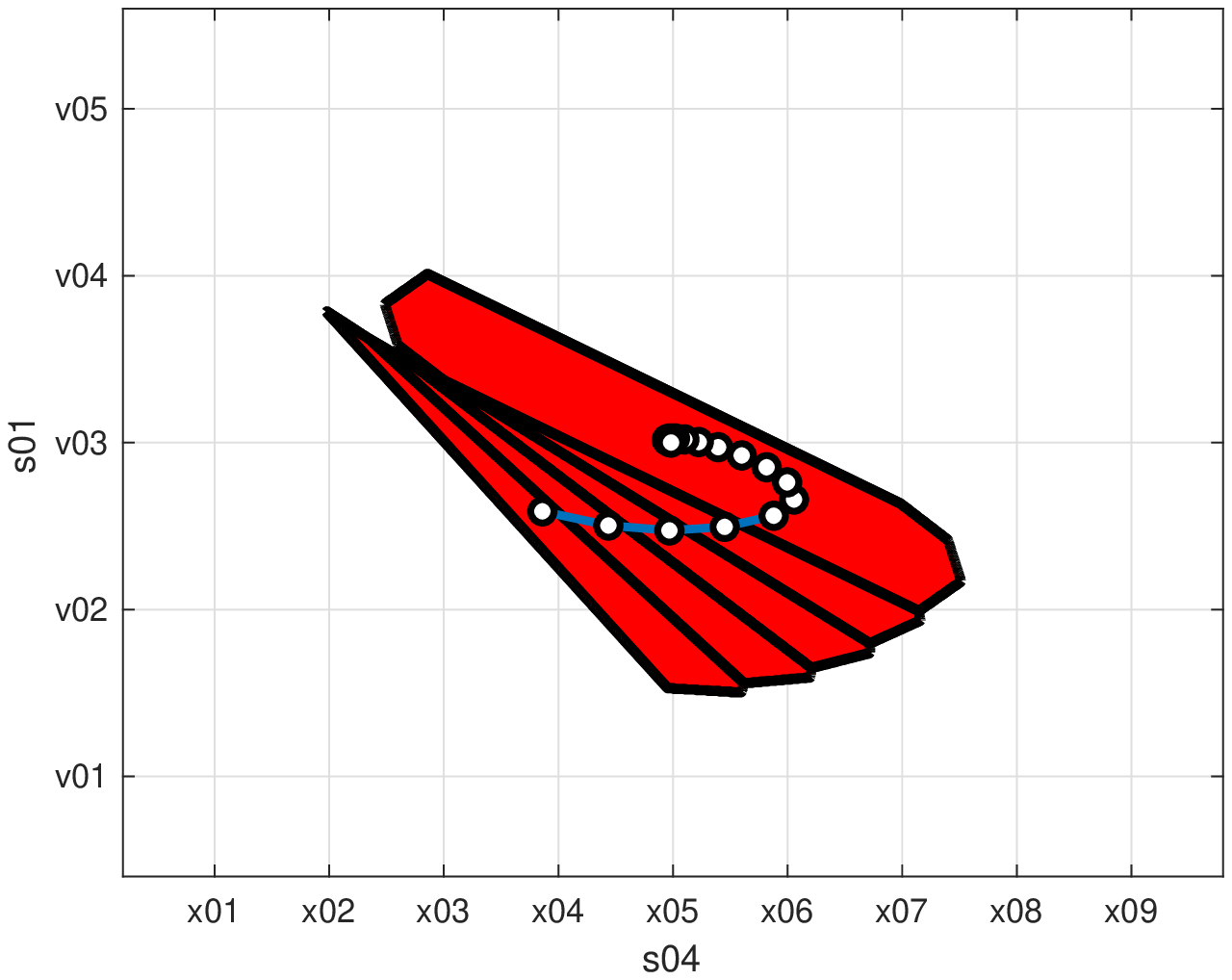}%
\end{psfrags}%
%
} \quad \quad \subfloat[Suboptimal approach]{
%
%
\begin{psfrags}%
\psfragscanon%
\scriptsize%
%
\psfrag{s01}[t][t]{\color[rgb]{0.15,0.15,0.15}\setlength{\tabcolsep}{0pt}\begin{tabular}{c}$x_1$\end{tabular}}%
\psfrag{s02}[b][b]{\color[rgb]{0.15,0.15,0.15}\setlength{\tabcolsep}{0pt}\begin{tabular}{c}$x_2$\end{tabular}}%
%
\color[rgb]{0.15,0.15,0.15}%
%
\psfrag{x01}[t][t]{-2}%
\psfrag{x02}[t][t]{-1.5}%
\psfrag{x03}[t][t]{-1}%
\psfrag{x04}[t][t]{-0.5}%
\psfrag{x05}[t][t]{0}%
\psfrag{x06}[t][t]{0.5}%
\psfrag{x07}[t][t]{1}%
\psfrag{x08}[t][t]{1.5}%
\psfrag{x09}[t][t]{2}%
%
\psfrag{v01}[r][r]{-1}%
\psfrag{v02}[r][r]{-0.5}%
\psfrag{v03}[r][r]{0}%
\psfrag{v04}[r][r]{0.5}%
\psfrag{v05}[r][r]{1}%
%
\includegraphics[width = 0.28\textwidth]{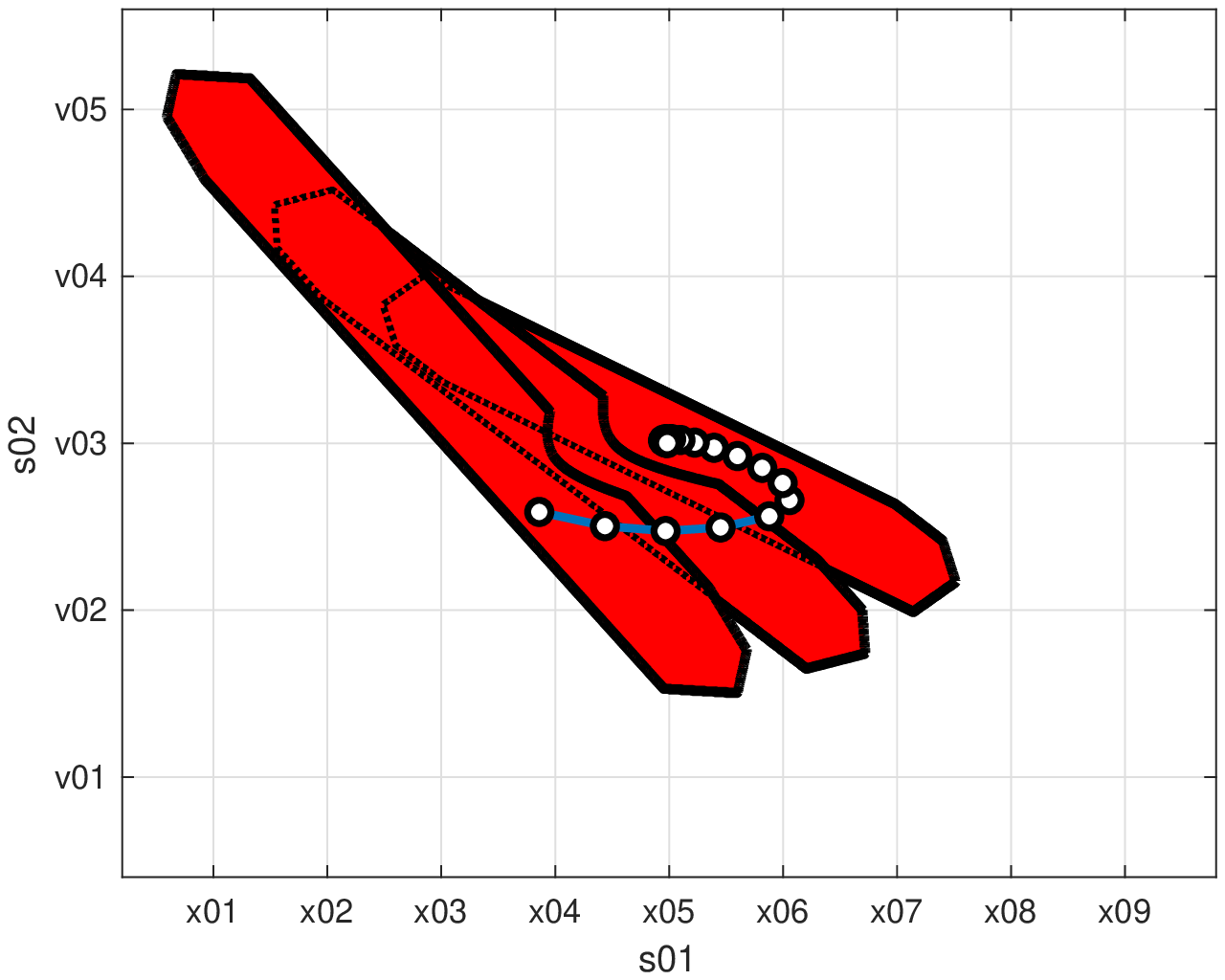}%
\end{psfrags}%
%
} \quad \quad \subfloat[Suboptimal approach (with projections)]{
%
%
\begin{psfrags}%
\psfragscanon%
\scriptsize%
%
\psfrag{s01}[t][t]{\color[rgb]{0.15,0.15,0.15}\setlength{\tabcolsep}{0pt}\begin{tabular}{c}$x_1$\end{tabular}}%
\psfrag{s03}[b][b]{\color[rgb]{0.15,0.15,0.15}\setlength{\tabcolsep}{0pt}\begin{tabular}{c}$x_2$\end{tabular}}%
%
\color[rgb]{0.15,0.15,0.15}%
%
\psfrag{x01}[t][t]{-2}%
\psfrag{x02}[t][t]{-1.5}%
\psfrag{x03}[t][t]{-1}%
\psfrag{x04}[t][t]{-0.5}%
\psfrag{x05}[t][t]{0}%
\psfrag{x06}[t][t]{0.5}%
\psfrag{x07}[t][t]{1}%
\psfrag{x08}[t][t]{1.5}%
\psfrag{x09}[t][t]{2}%
%
\psfrag{v01}[r][r]{-1}%
\psfrag{v02}[r][r]{-0.5}%
\psfrag{v03}[r][r]{0}%
\psfrag{v04}[r][r]{0.5}%
\psfrag{v05}[r][r]{1}%
%
\includegraphics[width = 0.28\textwidth]{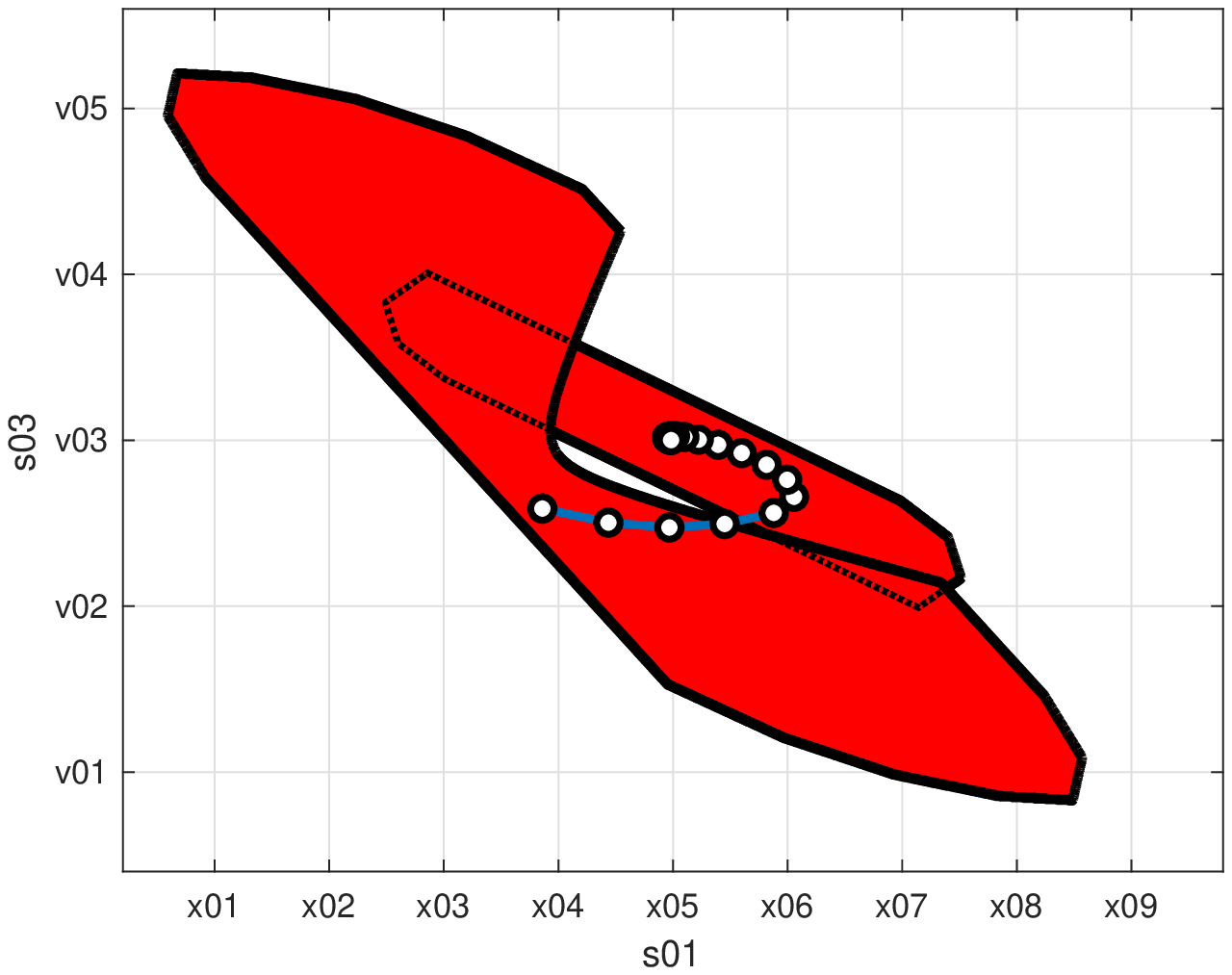}%
\end{psfrags}%
%
} 
\caption{Trajectories of the closed-loop system for a random initial state. The figure compares the optimal approach (left) with the new suboptimal approach (middle) and the new suboptimal approach using offline projections (right). } 
\label{fig:BeispielSISO}
\end{figure*}  

\begin{figure}[]
\centering
\captionsetup[subfloat]{justification=centering}
\subfloat[Optimal approach]{
%
%
\begin{psfrags}%
\psfragscanon%
\scriptsize%
%
\psfrag{s02}[b][b]{\color[rgb]{0.15,0.15,0.15}\setlength{\tabcolsep}{0pt}\begin{tabular}{c}$x(k)$\end{tabular}}%
\psfrag{s05}[b][b]{\color[rgb]{0.15,0.15,0.15}\setlength{\tabcolsep}{0pt}\begin{tabular}{c}$u(k)$\end{tabular}}%
\psfrag{s07}[t][t]{\color[rgb]{0.15,0.15,0.15}\setlength{\tabcolsep}{0pt}\begin{tabular}{c}$k$\end{tabular}}%
\psfrag{s08}[b][b]{\color[rgb]{0.15,0.15,0.15}\setlength{\tabcolsep}{0pt}\begin{tabular}{c}$e(k)$\end{tabular}}%
%
\color[rgb]{0.15,0.15,0.15}%
%
\psfrag{x01}[t][t]{0}%
\psfrag{x02}[t][t]{10}%
\psfrag{x03}[t][t]{20}%
\psfrag{x04}[t][t]{30}%
\psfrag{x05}[t][t]{40}%
\psfrag{x06}[t][t]{50}%
\psfrag{x07}[t][t]{0}%
\psfrag{x08}[t][t]{10}%
\psfrag{x09}[t][t]{20}%
\psfrag{x10}[t][t]{30}%
\psfrag{x11}[t][t]{40}%
\psfrag{x12}[t][t]{50}%
\psfrag{x13}[t][t]{0}%
\psfrag{x14}[t][t]{10}%
\psfrag{x15}[t][t]{20}%
\psfrag{x16}[t][t]{30}%
\psfrag{x17}[t][t]{40}%
\psfrag{x18}[t][t]{50}%
%
\psfrag{v01}[r][r]{0}%
\psfrag{v02}[r][r]{}%
\psfrag{v03}[r][r]{1}%
\psfrag{v04}[r][r]{-5}%
\psfrag{v05}[r][r]{0}%
\psfrag{v06}[r][r]{5}%
\psfrag{v07}[r][r]{-15}%
\psfrag{v08}[r][r]{0}%
\psfrag{v09}[r][r]{15}%
%
\includegraphics[width = 0.21\textwidth]{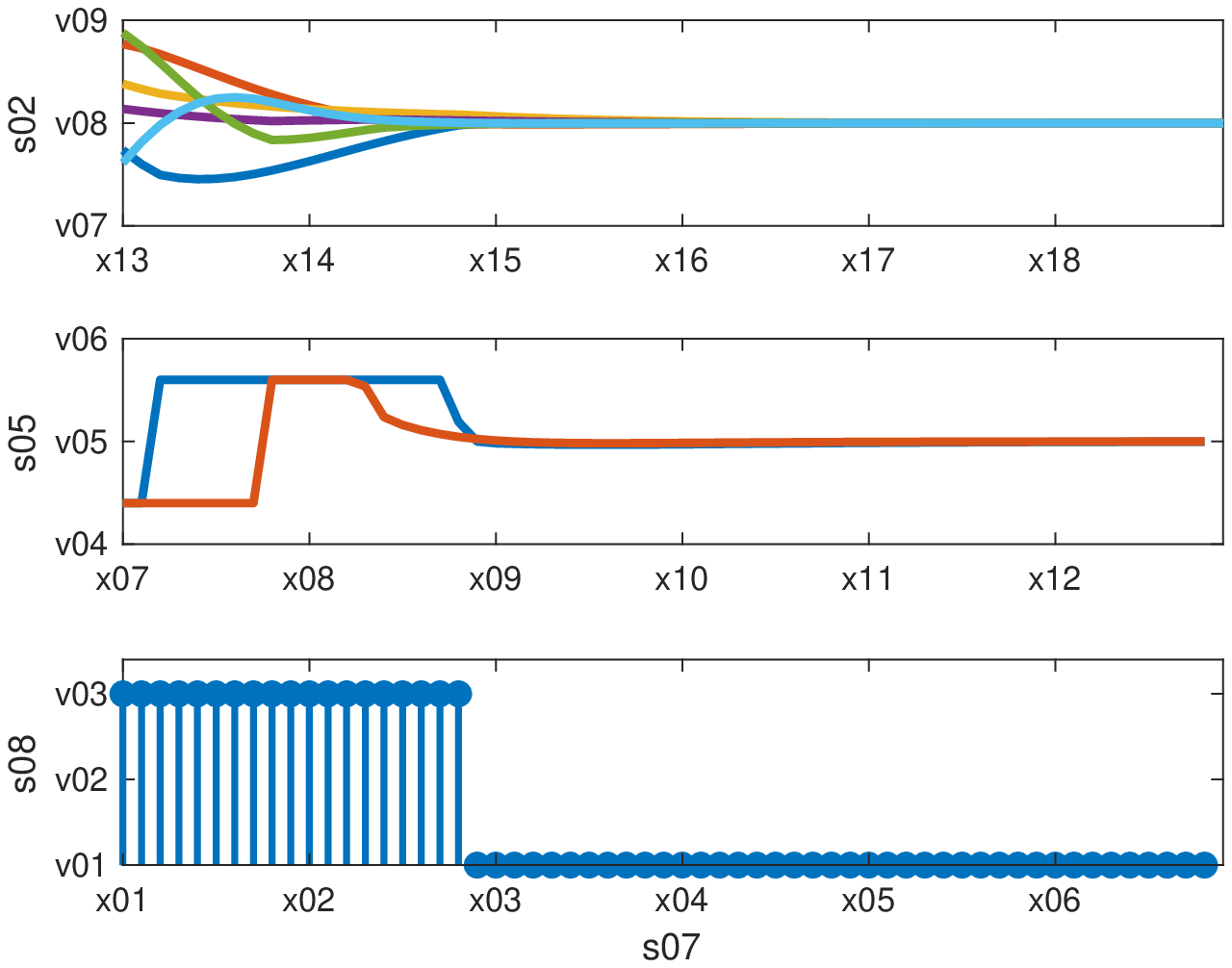}%
\end{psfrags}%
%
} \quad \quad
\subfloat[Suboptimal approach ($\lambda=1$)]{
%
%
\begin{psfrags}%
\psfragscanon%
\scriptsize%
%
\psfrag{s02}[b][b]{\color[rgb]{0.15,0.15,0.15}\setlength{\tabcolsep}{0pt}\begin{tabular}{c}$x(k)$\end{tabular}}%
\psfrag{s05}[b][b]{\color[rgb]{0.15,0.15,0.15}\setlength{\tabcolsep}{0pt}\begin{tabular}{c}$u(k)$\end{tabular}}%
\psfrag{s07}[t][t]{\color[rgb]{0.15,0.15,0.15}\setlength{\tabcolsep}{0pt}\begin{tabular}{c}$k$\end{tabular}}%
\psfrag{s08}[b][b]{\color[rgb]{0.15,0.15,0.15}\setlength{\tabcolsep}{0pt}\begin{tabular}{c}$e(k)$\end{tabular}}%
%
\color[rgb]{0.15,0.15,0.15}%
%
\psfrag{x01}[t][t]{0}%
\psfrag{x02}[t][t]{10}%
\psfrag{x03}[t][t]{20}%
\psfrag{x04}[t][t]{30}%
\psfrag{x05}[t][t]{40}%
\psfrag{x06}[t][t]{50}%
\psfrag{x07}[t][t]{0}%
\psfrag{x08}[t][t]{10}%
\psfrag{x09}[t][t]{20}%
\psfrag{x10}[t][t]{30}%
\psfrag{x11}[t][t]{40}%
\psfrag{x12}[t][t]{50}%
\psfrag{x13}[t][t]{0}%
\psfrag{x14}[t][t]{10}%
\psfrag{x15}[t][t]{20}%
\psfrag{x16}[t][t]{30}%
\psfrag{x17}[t][t]{40}%
\psfrag{x18}[t][t]{50}%
%
\psfrag{v01}[r][r]{0}%
\psfrag{v02}[r][r]{}%
\psfrag{v03}[r][r]{1}%
\psfrag{v04}[r][r]{-5}%
\psfrag{v05}[r][r]{0}%
\psfrag{v06}[r][r]{5}%
\psfrag{v07}[r][r]{-15}%
\psfrag{v08}[r][r]{0}%
\psfrag{v09}[r][r]{15}%
%
\includegraphics[width = 0.21\textwidth]{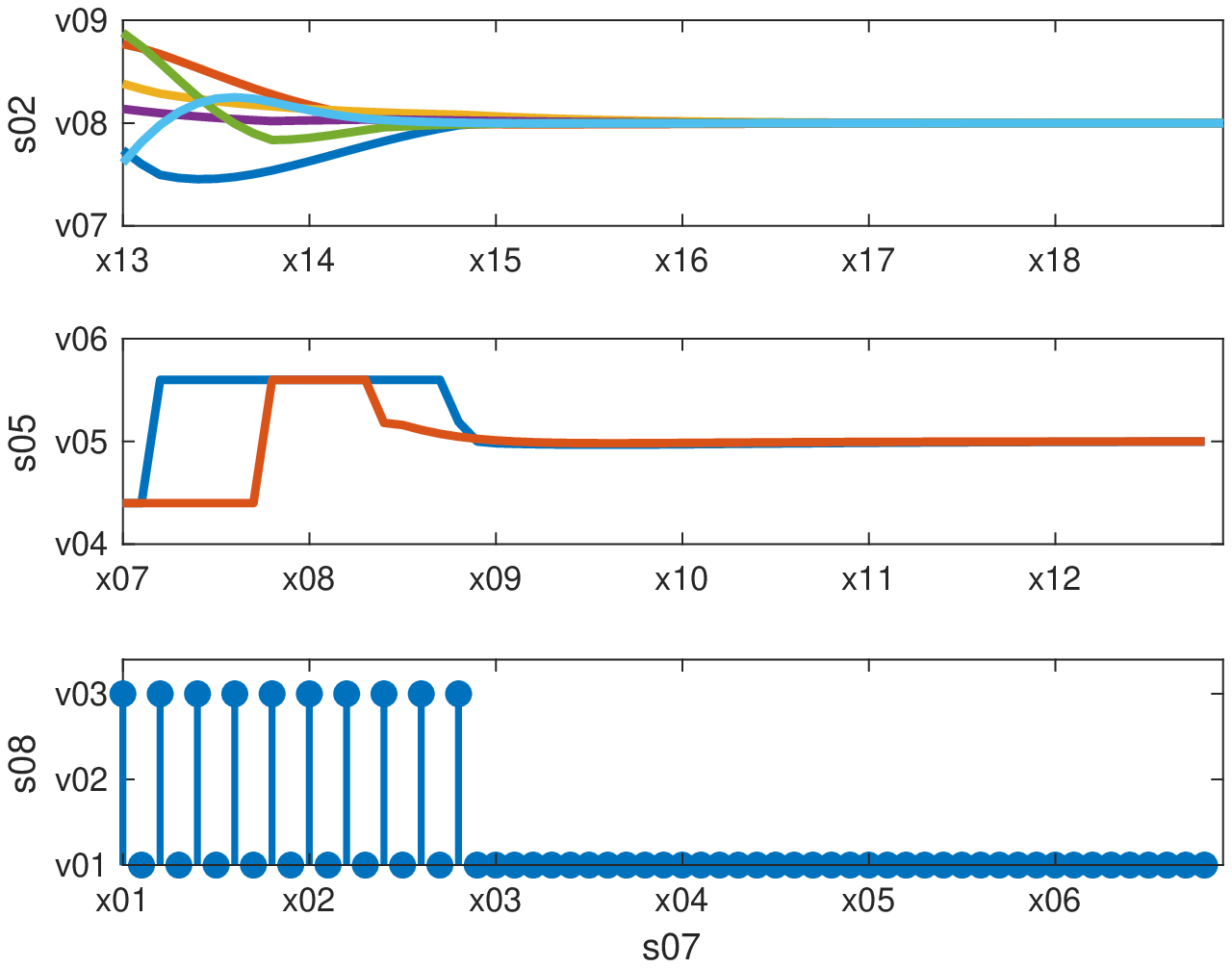}%
\end{psfrags}%
%
} \\ 
\subfloat[Suboptimal approach ($\lambda=1$, with projections)]{
%
%
\begin{psfrags}%
\psfragscanon%
\scriptsize%
%
\psfrag{s01}[t][t]{\color[rgb]{0.15,0.15,0.15}\setlength{\tabcolsep}{0pt}\begin{tabular}{c}$k$\end{tabular}}%
\psfrag{s02}[b][b]{\color[rgb]{0.15,0.15,0.15}\setlength{\tabcolsep}{0pt}\begin{tabular}{c}$e(k)$\end{tabular}}%
\psfrag{s05}[b][b]{\color[rgb]{0.15,0.15,0.15}\setlength{\tabcolsep}{0pt}\begin{tabular}{c}$x(k)$\end{tabular}}%
\psfrag{s11}[b][b]{\color[rgb]{0.15,0.15,0.15}\setlength{\tabcolsep}{0pt}\begin{tabular}{c}$u(k)$\end{tabular}}%
%
\color[rgb]{0.15,0.15,0.15}%
%
\psfrag{x01}[t][t]{0}%
\psfrag{x02}[t][t]{10}%
\psfrag{x03}[t][t]{20}%
\psfrag{x04}[t][t]{30}%
\psfrag{x05}[t][t]{40}%
\psfrag{x06}[t][t]{50}%
\psfrag{x07}[t][t]{60}%
\psfrag{x08}[t][t]{70}%
\psfrag{x09}[t][t]{80}%
\psfrag{x10}[t][t]{0}%
\psfrag{x11}[t][t]{10}%
\psfrag{x12}[t][t]{20}%
\psfrag{x13}[t][t]{30}%
\psfrag{x14}[t][t]{40}%
\psfrag{x15}[t][t]{50}%
\psfrag{x16}[t][t]{60}%
\psfrag{x17}[t][t]{70}%
\psfrag{x18}[t][t]{80}%
\psfrag{x19}[t][t]{0}%
\psfrag{x20}[t][t]{10}%
\psfrag{x21}[t][t]{20}%
\psfrag{x22}[t][t]{30}%
\psfrag{x23}[t][t]{40}%
\psfrag{x24}[t][t]{50}%
\psfrag{x25}[t][t]{60}%
\psfrag{x26}[t][t]{70}%
\psfrag{x27}[t][t]{80}%
%
\psfrag{v01}[r][r]{0}%
\psfrag{v02}[r][r]{}%
\psfrag{v03}[r][r]{1}%
\psfrag{v04}[r][r]{-5}%
\psfrag{v05}[r][r]{0}%
\psfrag{v06}[r][r]{5}%
\psfrag{v07}[r][r]{-15}%
\psfrag{v08}[r][r]{0}%
\psfrag{v09}[r][r]{15}%
%
\includegraphics[width = 0.21\textwidth]{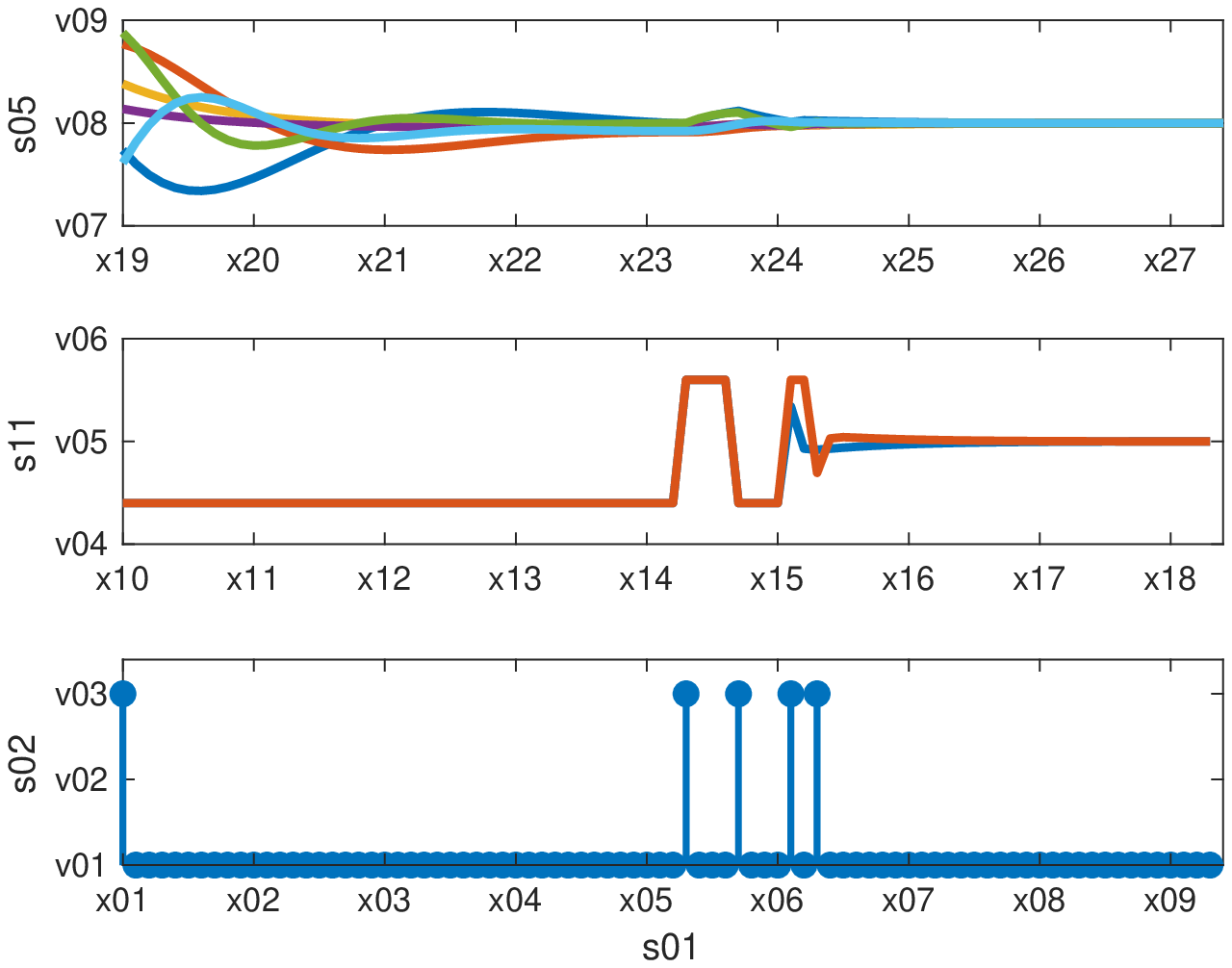}%
\end{psfrags}%
%
} \quad \quad
\subfloat[Suboptimal approach ($\lambda=0.8$, with projections)]{
%
%
\begin{psfrags}%
\psfragscanon%
\scriptsize%
%
\psfrag{s02}[b][b]{\color[rgb]{0.15,0.15,0.15}\setlength{\tabcolsep}{0pt}\begin{tabular}{c}$x(k)$\end{tabular}}%
\psfrag{s05}[b][b]{\color[rgb]{0.15,0.15,0.15}\setlength{\tabcolsep}{0pt}\begin{tabular}{c}$u(k)$\end{tabular}}%
\psfrag{s07}[t][t]{\color[rgb]{0.15,0.15,0.15}\setlength{\tabcolsep}{0pt}\begin{tabular}{c}$k$\end{tabular}}%
\psfrag{s08}[b][b]{\color[rgb]{0.15,0.15,0.15}\setlength{\tabcolsep}{0pt}\begin{tabular}{c}$e(k)$\end{tabular}}%
%
\color[rgb]{0.15,0.15,0.15}%
%
\psfrag{x01}[t][t]{0}%
\psfrag{x02}[t][t]{10}%
\psfrag{x03}[t][t]{20}%
\psfrag{x04}[t][t]{30}%
\psfrag{x05}[t][t]{40}%
\psfrag{x06}[t][t]{50}%
\psfrag{x07}[t][t]{0}%
\psfrag{x08}[t][t]{10}%
\psfrag{x09}[t][t]{20}%
\psfrag{x10}[t][t]{30}%
\psfrag{x11}[t][t]{40}%
\psfrag{x12}[t][t]{50}%
\psfrag{x13}[t][t]{0}%
\psfrag{x14}[t][t]{10}%
\psfrag{x15}[t][t]{20}%
\psfrag{x16}[t][t]{30}%
\psfrag{x17}[t][t]{40}%
\psfrag{x18}[t][t]{50}%
%
\psfrag{v01}[r][r]{0}%
\psfrag{v02}[r][r]{}%
\psfrag{v03}[r][r]{1}%
\psfrag{v04}[r][r]{-5}%
\psfrag{v05}[r][r]{0}%
\psfrag{v06}[r][r]{5}%
\psfrag{v07}[r][r]{-15}%
\psfrag{v08}[r][r]{0}%
\psfrag{v09}[r][r]{15}%
%
\includegraphics[width = 0.21\textwidth]{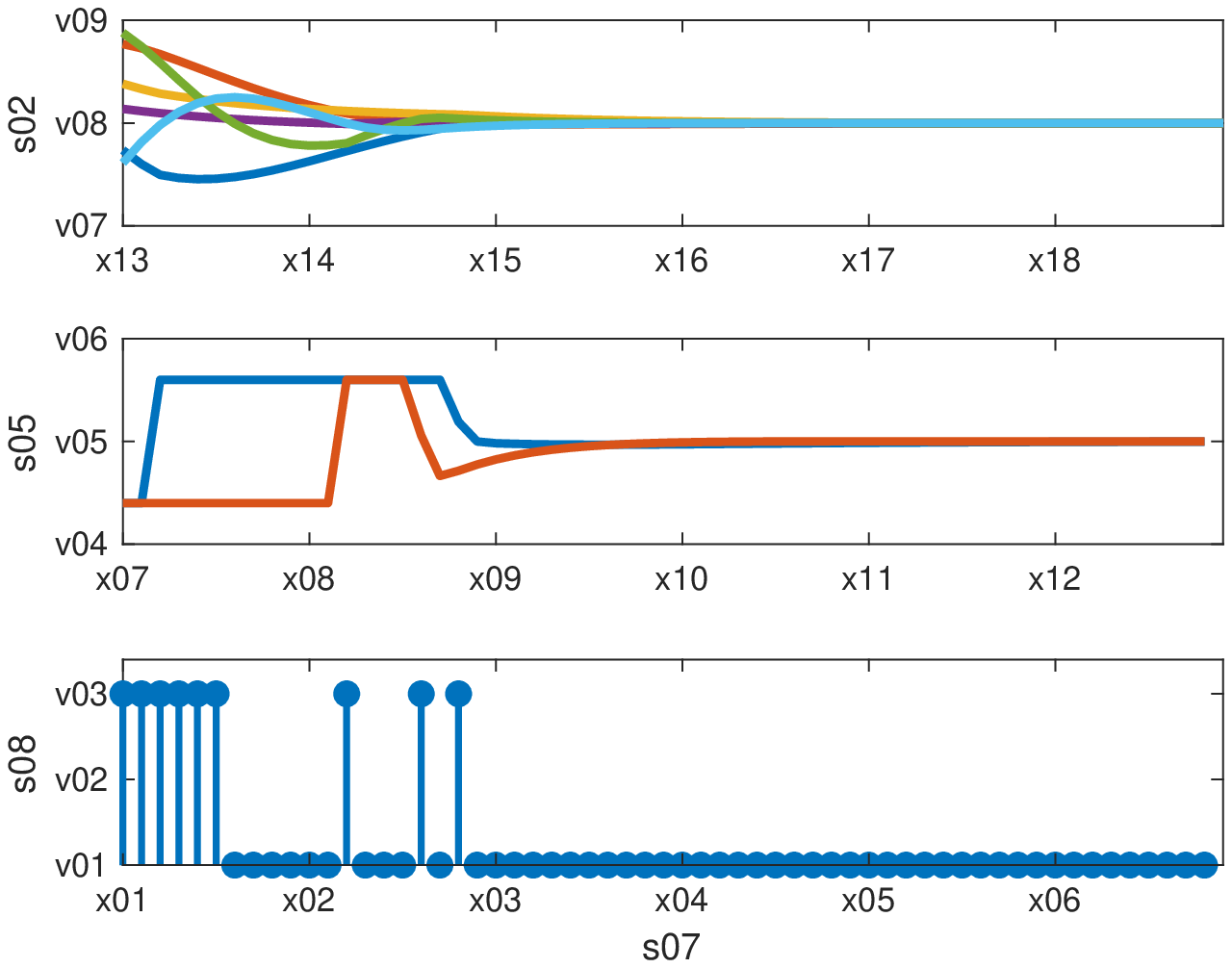}%
\end{psfrags}%
%
} \\
\caption{Trajectories of the closed-loop system for a random initial state. The figure compares the optimal approach (a) with the new suboptimal approach (b) and the new suboptimal approach using offline projections (c, d). }  
\label{fig:BeispielMIMO}
\end{figure}  

Fig. \ref{fig:BeispielSISO} shows the results of the closed-loop system for a random initial state for Example 1. We compare the original regional MPC approach (left) to our new approach (middle and right). The results depicted in the middle and right plot differ with respect to the feasibility region that is used. The middle plot uses feasibility regions $\mathcal{B}=\mathcal{F}$ according to Proposition \ref{prop:Feasibility}. 
The right plot additionally uses feasibility regions  $\mathcal{B}=\mathcal{C}$ according to Proposition \ref{prop:Cl}. Feasibility regions $\mathcal{C}$ result by projection. Since online projections are computationally demanding, we here only use projections for saturated feedback laws that exist in the solution to \eqref{eq:ReformulatedMPCProblem}. Note that these laws can be easily found a priori such that the corresponding projections can be computed offline before.    

The figure shows the trajectories of the states $x(k)$, the inputs $u(k)$ and the indicator function $e(k)$ with $e(k)=1$ if a QP is solved and $e(k)=0$ otherwise (upper plots). It also illustrates the validity regions in state-space (bottom plots). \\ 
It is obvious from the figure that each of the approaches regulates the system to the origin. 
In the original approach, five QPs and thus five feedback laws and polytopes must be calculated along the closed-loop trajectory. Only the feedback law associated with the last polytope can be reused for more than one time step. The nonlinearly bounded regions of our new approach are depicted in the middle and right plot. The figures show that the new regions are larger than the optimal polytopes. 
With the suboptimal approach (middle) only three feedback laws and thus three regions must be calculated to regulate the system to the origin. The trajectories of the states and inputs are very similar to the trajectories in the optimal case. With the suboptimal approach using projections (right), only two QPs and thus two regions must be calculated and the trajectories are again very similar to the optimal case. The first region is larger than the first region calculated in the suboptimal approach without projections since the feasibility region is larger. \\
Since results on any single trajectory are only anecdotal, Table \ref{tab:Ergebnisse} states the results for 1000 random initial states. Again the optimal approach is compared to the new suboptimal approach. The second column shows the results for the optimal approach. The results for the new approach without offline projections is stated in the third column. The fourth and fifth column show the results of the new approach extended by offline projections for saturated feedback laws. These results differ with respect to the choice of the parameter $\lambda$. 
The table states the number of saved QPs ($\Delta$QP) on the central node, the number of saved flops ($\Delta$flops) on the local node and the deviation of the total costs ($\Delta$costs) in percent compared to the optimal approach along the closed-loop trajectories. The costs are calculated by summing up the cost function in \eqref{eq:ReformulatedMPCProblem} in every time step along the closed-loop trajectories. They give information about the deviation of the trajectories of the suboptimal approach from the optimal approach and are utilized to measure the performance of the approaches.\\
In the suboptimal approach without projections and $\lambda=1$, the number of solved QPs on the central node is reduced by about \unit[24.71]{\%} and the number of flops on the local node is increased by about \unit[9.37]{\%}.The costs are increased only by about \unit[1.08]{\%}. With the suboptimal approach using projections and $\lambda=1$ the number of solved QPs can even be reduced further, specifically by about \unit[33.48]{\%}. The number of flops is similar to the optimal approach and can even be reduced by about \unit[2.12]{\%}. 
However, the reduction in the number of QPs results in a worse performance compared to the suboptimal approach without projections resulting in an increase of the costs by about \unit[4.23]{\%}. By choosing $\lambda=0.8$ the performance can be improved compared to the case $\lambda=1$ but more QPs have to be solved in this case. Nevertheless, compared to the optimal case, about \unit[29.46]{\%} of the QPs are saved and the number of flops can be reduced by about \unit[2.4]{\%}. \\          

Figure \ref{fig:BeispielMIMO} shows the results of the closed-loop system for a random initial state for Example 2. Again we compare the optimal approach (a) to the suboptimal approach (b,c,d). Figure (b) illustrates the results for the suboptimal approach without projections and $\lambda=1$. In Figures (c) and (d) projections with $\lambda=1$ and $\lambda=0.8$, respectively, are used. As in the previous example, all approaches regulate the system to the origin. In the optimal approach (a) 19 QPs must be solved. In the suboptimal approach without projections and $\lambda=1$, i.e. in (b), the number of solved QPs can be reduced by 9 QPs compared to the optimal approach, while the trajectories are very similar to the optimal case. 
The suboptimal approach with projections and $\lambda=1$ (c) saves 14 QPs compared to the optimal approach. However,  although it has the largest savings potential, the trajectories deviate at most from the optimal trajectories in this case. Specifically, the terminal region is reached after only 53 time steps, in contrast to the 18 time steps required in the optimal case. This problem can be handled by decreasing the parameter $\lambda$. The suboptimal approach with projections and $\lambda=0.8$ results in trajectories that are similar to the trajectories in the optimal approach, while reducing the number of QPs by 10 QPs. 
Again Table \ref{tab:Ergebnisse} states the results for 1000 random initial states. The suboptimal approach without projections and $\lambda=1$ saves about \unit[23.67]{\%} of the QPs on the central node. The number of flops on the local node is increased by only \unit[1.05]{\%}.The costs are increased only by about \unit[0.02]{\%}. With the suboptimal approach using projections and $\lambda=1$ the number of QPs can be reduced by about \unit[73.52]{\%}. Moreover, the flops on the local node can be reduced by about \unit[37.21]{\%}. 
However, the reduction in the number of QPs results in a large deviation from the optimal trajectories by about \unit[55.79]{\%}. By choosing $\lambda=0.8$ the performance can be improved  compared to the case $\lambda=1$ but more QPs have to be solved in this case. Nevertheless, compared to the optimal case about \unit[50.63]{\%} of the QPs are skipped and the number of flops can be reduced by about \unit[23.78]{\%}. \\     

\setlength{\tabcolsep}{2pt}
\begin{table*}[t!]
\caption{Computational effort of the new approach in a networked control setting. Results are shown for 1000 random initial states. The new approach is compared to the approach from \cite{Jost2015a}. The table states the number of saved QPs on the central node ($\Delta$QPs), the number of saved flops on the local node ($\Delta$flops) and the deviation of the costs ($\Delta$costs).}
\begin{center}
\resizebox{15.5cm}{!}{%
\begin{tabular}{|c|ccc|ccc|ccc|ccc|}
\hline & \multicolumn{3}{|c|}{\multirow{2}{*}{Optimal Approach}} & \multicolumn{9}{|c|}{Suboptimal Approach} \\  \cline{5-13}
 Sys. & \multicolumn{3}{|c|}{} & \multicolumn{3}{|c|}{$(\lambda=1, \text{no projections})$}&  \multicolumn{3}{|c|}{$(\lambda=1 \ \text{with projections})$} &  \multicolumn{3}{|c|}{$(\lambda=0.8 \ \text{with projections})$} \\ \cline{2-13}
   & QPs & flops & costs & $\Delta$QPs & $\Delta$flops & $\Delta$costs & $\Delta$QPs &  $\Delta$flops & $\Delta$costs  & $\Delta$QPs &  $\Delta$flops & $\Delta$costs \\  \hline
Ex. 1 & 4471 & 2005632 & 4932.35  &   -\unit{24.71}{\%}  &  +\unit{9.37}{\%}   & \unit{1.08}{\%} & -\unit{33.48}{\%} & -\unit{2.12}{\%} & \unit{4.23}{\%} & -\unit{29.46}{\%} & -\unit{2.4}{\%} & \unit{2.03}{\%} \\  \hline
Ex. 2 & 17910 & 454,93 $\cdot 10^6$  & 5,75 $\cdot 10^7$  & -\unit{23.67}{\%}  & +\unit{1.05}{\%}  &\unit{0.02}{\%} & -\unit{73.52}{\%}  & -\unit{37.21}{\%} & \unit{55.79}{\%}  & -\unit{50.63}{\%} & -\unit{23.78}{\%}  & \unit{1.02}{\%}  \\  \hline
\end{tabular}%
}
\label{tab:Ergebnisse}
\end{center}
\end{table*}

\section{Conclusions}\label{sec:conclusion}
We improved an event-triggered MPC approach, in which feedback laws defined on polytopic regions are calculated on a central node by solving a QP and evaluated on a local node as long as they provide the optimal solution. By dropping the requirement for optimality, we calculated an extended region of validity for a feedback law that is simple to evaluate in a networked control setting. Feasibility and stability are maintained and these requirements result in a nonlinearly bounded region defined by simple linear and quadratic inequalities. 
We showed that the computational effort on the central and local node can be reduced compared to the existing approach that enforces optimality. 
Furthermore, we showed that the new approach can be adjusted for a desired closed-loop performance with a easy-to-interpret tuning parameter.       

\section*{Acknowledgement}                           
Support by the Deutsche Forschungsgemeinschaft (DFG) under grant MO 1086/15-1 is gratefully acknowledged.

\bibliographystyle{plain}      
\bibliography{literature.bib}

  \small               
                                
\end{document}